%
%

\magnification=1200

\font\titfont=cmr10 at 12 pt

\font \fr = eufm10


\font\AAA=cmr14 at 12pt
\font\BBB=cmr14 at 8pt

\overfullrule=0in

\def\oa#1{\overrightarrow #1}
\def\dim{{\rm dim}}
\def\dist{{\rm dist}}

\def\deg{{\rm deg}}

\def\log{{\rm log}}
\def\Hess{{\rm Hess}}
\def\Sym{\rm Sym}

\def\tr{{\rm tr}}
\def\max{{\rm max}}

\def\span{{\rm span\,}}
\def\Hom{{\rm Hom\,}}

\def\End{{\rm End}}

\def\arr{\longrightarrow}
\def\supp{{\rm supp}}

\def\foral{\qquad {\rm for\ all\ \ }}
\def\fpsh{{\cal PSH}(X,\f)}
\def\Core{{\rm Core}}
\def\dis{f_M}


\def\Theorem#1{\medskip\noindent {\AAA T\BBB HEOREM \rm #1.}}
\def\Prop#1{\medskip\noindent {\AAA P\BBB ROPOSITION \rm  #1.}}
\def\Cor#1{\medskip\noindent {\AAA C\BBB OROLLARY \rm #1.}}
\def\Lemma#1{\medskip\noindent {\AAA L\BBB EMMA \rm  #1.}}
\def\Remark#1{\medskip\noindent {\AAA R\BBB EMARK \rm  #1.}}
\def\Note#1{\medskip\noindent {\AAA N\BBB OTE \rm  #1.}}
\def\Def#1{\medskip\noindent {\AAA D\BBB EFINITION \rm  #1.}}

\def\Ex#1{\medskip\noindent {\AAA E\BBB XAMPLE \rm    #1.}}
\def\Qu#1{\medskip\noindent {\AAA Q\BBB UESTION \rm    #1.}}

\def\HoQu#1{ {\AAA T\BBB HE\ \AAA H\BBB ODGE\ \AAA Q\BBB UESTION \rm    #1.}}

\def\pf{\medskip\noindent {\bf Proof.}\ }
\def\qed{\hfill  $\vrule width5pt height5pt depth0pt$}

\def\qedqed{\hfill  $\vrule width5pt height5pt depth0pt$ $\vrule width5pt height5pt depth0pt$}

\def\df{d^{\phi}}
\def\hk{\_{\rm l}\,}
\def\n{\nabla}
\def\w{\wedge}

   \def\cc{{\cal C}}     
   \def\cp{{\cal P}}
   
\def\ce{{\cal E}}   
\def\ch{{\cal H}}   \def\cm{{\cal M}}
\def\cs{{\cal S}}   \def\cn{{\cal N}}
\def\cd{{\cal D}}
\def\cl{{\cal L}}
\def\cp{{\cal P}}
\def\cf{{\cal F}}
\def\ccr{{\cal  R}}

\def\gerG{{\fr{\hbox{g}}}}

\def\vf{\varphi}

\def\wt{\widetilde}
\def\wh{\widehat}

\def\and{\qquad {\rm and} \qquad}
\def\arr{\longrightarrow}
\def\ol{\overline}
\def\bbr{{\bf R}}\def\bbh{{\bf H}}\def\bbo{{\bf O}}
\def\bbc{{\bf C}}

\def\bbz{{\bf Z}}

\def\a{\alpha}
\def\b{\beta}
\def\d{\delta}
\def\e{\epsilon}
\def\f{\phi}

\def\o{\omega}

\def\s{\sigma}
\def\x{\xi}
\def\z{\zeta}

\def\D{\Delta}
\def\L{\Lambda}
\def\G{\Gamma}
\def\O{\Omega}

\def\bd{\partial}
\def\bdf{\partial_{\f}}
\def\fp{$\phi$-plurisubharmonic }
\def\fh{$\phi$-pluriharmonic }

\def\ffl{$\f$-flat}
\def\PH#1{\widehat {#1}}

\def\lp{\Lambda_+(\f)}
\def\lpp{\Lambda^+(\f)}
\def\bo{\partial \Omega}
\def\fc{$\phi$-convex }
\def\PSH{ \cp\cs\ch}
\def\totr{ $\phi$-free }
\def\BM{\lambda}
\def\Der{D}
\def\CH{{\cal H}}
\def\RH{\overline{\ch}^\f }

\def\B{A}
\def\YY{0}
\def\AA{1}
\def\ZZ{2}
\def\BB{3}
\def\CC{4}
\def\DD{5}
\def\EE{6}
\def\FF{7}
\def\GG{8}
\def\HH{9}
\def\II{10}

\ 
\vskip .3in

\centerline{\titfont PLURISUBHARMONIC FUNCTIONS IN }
\smallskip

\centerline{\titfont CALIBRATED GEOMETRIES }
\bigskip

\centerline{\titfont F. Reese Harvey and H. Blaine Lawson, Jr.$^*$}
\vglue .9cm
\smallbreak\footnote{}{ $ {} \sp{ *}{\rm Partially}$  supported by
the N.S.F. }

\vskip .3in
\centerline{\bf ABSTRACT} \medskip
  \font\abstractfont=cmr10 at 10 pt

{{\parindent= .7in\narrower\abstractfont \noindent
In this paper we introduce and study the notion of plurisubharmonic functions
in  calibrated geometry.  These functions generalize 
the classical plurisubharmonic functions from complex geometry
 and enjoy many of their important properties. 
 Moreover, they exist in abundance whereas the corresponding
 pluriharmonics are generally quite scarce.
 A number of the results
established in complex analysis via plurisubharmonic functions are extended
to calibrated manifolds.  This paper investigates, in depth, questions
of:  pseudo-convexity and cores,  positive $\f$-currents, 
Duval-Sibony Duality, and boundaries of $\f$-submanifolds,
all in the context of a general calibrated manifold $(X,\f)$.
Analogues of totally real submanifolds are used to construct
enormous families of strictly $\f$-convex
spaces  with every topological type allowed by 
Morse Theory. Specific calibrations are used as examples throughout.
Analogues of the Hodge Conjecture in  calibrated geometry are considered.

}}

\vfill\eject\
\vskip 1in

\centerline{\bf TABLE OF CONTENTS} \bigskip

\qquad \YY. Introduction.\smallskip

\qquad \AA. Plurisubharmonic Functions.

\qquad\qquad\qquad\qquad Appendix: Pluriharmonic Functions.
\smallskip

\qquad \ZZ.  The $\f$-Hessian.  

\qquad\qquad\qquad\qquad  Appendix A: Submanifolds which are $\f$-Critical.

\qquad\qquad\qquad\qquad  Appendix B: Constructing $\f$-Plurisubharmonic Functions.
\smallskip

\qquad \BB.  Convexity in Calibrated Geometries.

\qquad\qquad \qquad\qquad Appendix A: Structure of the Core.

\qquad\qquad \qquad\qquad Appendix B: Examples of Complete Convex Manifolds and Cores.
\smallskip

\qquad \CC. Boundary Convexity.\smallskip

\qquad \DD.  Positive Currents in Calibrated Geometries. 

\qquad\qquad\qquad\qquad Appendix:  {The Reduced $\f$-Hessian. }  \smallskip 

\qquad \EE.  Duval-Sibony Duality.\smallskip

\qquad \FF. $\f$-Free Submanifolds..\smallskip

\qquad \GG.  Hodge Manifolds.\smallskip

\qquad \HH. Boundaries of $\f$-Submanifolds.\smallskip

\qquad \II. $\f$-Flat Hypersurfaces and Functions which are $\f$-Pluriharmonic mod $d$.\smallskip

\vfill\eject

\centerline{\bf \YY. Introduction.}
\medskip

Calibrated geometries, as introduced in [HL$_3$],  are  geometries
of distinguished submanifolds  determined by a fixed, closed differential form $\f$ on
a riemannian manifold $X$.
The basic example is that of a K\"ahler manifold (or more generally a symplectic manifold,
with compatible complex structure) where the distinguished submanifolds are the  holomorphic curves.
However, there exist many other interesting  geometries, each carrying a wealth of $\f$-submanifolds, particularly  on spaces with special holonomy.  These have attracted particular attention in recent
years due to their appearance in generalized Donaldson theories and in modern versions of string theory in Physics.

Unfortunately, analysis on these spaces $(X,\f)$ has been difficult,  in part because  
there is generally no reasonable analogue of
the  holomorphic functions  and  transformations which exist in the  K\"ahler  case.
However, in complex analysis there are many important results which can be established
using only the plurisubharmonic  functions. It turns out that analogues of these functions 
exist in abundance on any calibrated manifold, and they enjoy almost all the pleasant properties 
of their cousins from complex analysis. The point of this paper is to introduce and study 
these functions.

We begin by defining our notion of a  $\f$-plurisubharmonic function
on any calibrated manifold $(X,\f)$.  In the K\"ahler case  they
are exactly the classical plurisubharmonic functions.  We then study the basic properties of these functions, and subsequently use them to establish a series of results in geometry 
and analysis on $(X,\f)$.  

A  fundamental result is that:\smallskip

\centerline{\sl The restriction of a \fp function to a $\f$-submanifold $M$ is subharmonic }\smallskip

\noindent
in the induced metric on $M$.

Any  convex function on the riemannnian manifold $X$ is  $\f$-plurisubharmonic. Moreover, at least locally, there exists an abundance 
of \fp  functions which are not convex. 

The definition of  $\f$-plurisubharmonicity extends to arbitrary distributions on $X$.
For most calibrations -- including all the ``classical'' ones -- such distributions enjoy
all the nice properties of generalized subharmonic functions such as being  
locally Lebesgue integrable and represented by an upper-semicontinuous function taking
values in   $[-\infty, \infty)$. In general the maximum $F=\max\{f,g\}$ of two smooth \fp functions
is again \fp and can be  uniformly approximated by a decreasing sequence of smooth 
\fp functions.

To define \fp functions on a calibrated manifold $(X,\f)$ we introduce a second order
differential operator
$
\ch^\f : C^\infty(X) \ \to\ \ce^p(X),
$
the {\sl $\f$-Hessian}, given by  
$$
\ch^\f(f)\ =\ \BM_\f(\Hess f)
$$
where $\Hess f$ is the riemannian hessian of $f$ and $\BM_\f:\End(TX) \to \L^pT^*X$
is the bundle map given by  $\BM_\f(A) = D_{A^*}(\f)$ where  $D_{A^*}:\L^pT^*X\arr \L^pT^*X$ is the natural extension of $A^*:T^*X \to T^*X$ as a derivation.

When the calibration $\f$ is parallel there is a natural factorization
$$
\ch^\f \ =\ d d^\f
$$
where $d$ is the de Rham differential and $d^\f: C^\infty(X) \ \to\ \ce^{p-1}(X)$ is given by
$$
d^\f f\ \equiv \ \nabla f \hk \f.
$$
In general these operators are related by the equation: $\ch^\f  f = dd^\f f - \nabla_{\nabla f} (\f)$.

Recall that a calibration $\f$ of degree $p$  is a closed $p$-form  with the property that $\f(\x)\leq 1$ for all unit simple tangent $p$-vectors $\x$ on $X$. Those $\x$ for which $\f(\x)=1$ are called 
$\f$-planes, and the set of $\f$-planes is denoted by $G(\f)$. With this understood, a function
$f\in C^\infty(X)$ is defined to be {\bf  \fp} if $\ch^\f(f)(\x) \geq 0$ for all $\x\in G(\f)$. It is {\bf strictly
\fp } at a point $x\in X$ if $\ch^\f(f)(\x) > 0$ for all $\f$-planes $\x$ at $x$.  In a similar fashion, $f$ 
is called {\bf $\f$-pluriharmonic} if $\ch^\f(f)(\x) = 0$ for all $\x\in G(\f)$.
Denote by $\PSH(X,\f)$ the convex cone of \fp functions on $X$.

When $X$ is a complex manifold with a K\"ahler form $\o$, one easily computes that
$d^\o = d^c$, the conjugate differential.  In this case, $\ch^\o = dd^\o =dd^c$ and the 
$\o$-planes correspond to the complex lines in $TX$. Hence, the definitions above coincide
with the classical notions of plurisubharmonic and pluriharmonic functions on $X$.

With this said, we must remark that in many calibrated manifolds the $\f$-pluriharmonic
functions are scarce. For the  calibrations on manifolds with strict G$_2$ or Spin$_7$ holonomy,  for example,  every pluriharmonic function is constant.  For the Special Lagrangian
calibration $\f = {\rm Re}\{dz\}$, every $\f$-pluriharmonic  function $f$ defined locally in $\bbc^n$ is of the  form $f=a+q$ where $a$ is affine and $q$ is the real part of a complex quadratic form
(cf.  [Fu].) Nevertheless, as we stated above, the \fp functions in any calibrated geometry are locally abundant.

The fundamental property of the $\f$-Hessian: 
$$
\left(\ch^\f f\right)(\x)\ =\ {\rm trace}\left\{ \Hess f\bigr|_{\x}  \right\} \qquad \ \ {\rm for \ all \ } \f-{\rm planes\ } \x
$$
is established in Section \ZZ \ (Corollary \ZZ.5).

Beginning with Section \BB\  the \fp functions are used to study geometry and analysis 
on calibrated manifolds.  The first concept to be addressed is the analogue
of pseudoconvexity in complex geometry.

\medskip
\centerline
{\AAA C\BBB ONVEXITY.}\smallskip

 Let  $(X,\f)$ be a calibrated manifold and $K\subset X$ a closed
subset. By the {\bf $\f$-convex hull of $K$ } we mean the subset
$$
\wh K\ =\ \{x\in X : f(x) \leq \sup_K f \ \ {\rm for \ all \  } f\in \PSH(X,\f)\}
$$
The manifold $(X,\f)$ is said to be {\bf  $\f$-convex} if $K\subset \subset X \ \Rightarrow
\  \wh K \subset \subset X$  for all $K$.

\Theorem {\BB.3} {\sl A calibrated manifold  $(X,\f)$ is $\f$-convex if and only if 
it admits a \fp proper exhaustion function $f:X\to \bbr$.} \medskip

The manifold $(X,\f)$ will be called {\bf strictly $\f$-convex} if it admits an exhaustion function $f$ which is  strictly $\f$-plurisubharmonic, 
and it will be called  {\bf strictly $\f$-convex at infinity}  if $f$ is strictly \fp  outside of a compact subset.
It is shown that in the second case,  $f$ can be assumed to be
\fp everywhere. 
Analogues of Theorem \BB.3 are established in each of these cases.

Note that in complex geometry, \fc manifolds are Stein and manifolds which are \fc at infinity are 
 strongly pseudoconvex.

We next consider the {\bf core} of $X$ which is defined to be the set of points
$x\in X$ with the property that no $f\in\PSH(X,\f)$ is strictly \fp at $x$. The following results are established:
\smallskip

1)\ Every compact $\f$-submanifold is contained in Core$(X)$.   In  fact, every 
$\partial$-closed $\f$-positive current $T$ is supported in the core. More generally,
every $\f$-positive current $T$ is supported in $\wh {\supp(d T)} \cup \Core(X)$.
(See [HL$_3$] or \S\S \DD \ and \EE\  below for a discussion of $\f$-positive currents.)\smallskip

2)\ The manifold $X$ is strictly $\f$-convex at infinity if and only if Core$(X)$ is compact.
\smallskip

3)\ The manifold $X$ is strictly $\f$-convex if and only if Core$(X)=\emptyset$.
\smallskip

Examples of complete calibrated manifolds with compact  cores are given in an appendix 
to \S \BB. 
A very general construction of strictly \fc manifolds is presented in \S \GG.

We next examine the analogues of pseudoconvex boundaries in calibrated geometry.
\medskip
\centerline
{\AAA B\BBB OUNDARY \AAA C\BBB ONVEXITY.}\smallskip

Let $\O\subset X$ be a domain with smooth boundary $\bo$, and let $\rho:X\to \bbr$ be a {\bf defining function for $\partial \O$}, that is, a smooth function defined on a neighborhood of $\overline \O$ with $\O=\{x : \rho(x) <0\}$, and $\nabla \rho \neq 0$ on $\bo$. Then $\bo$ is said
to be {\bf $\f$-convex} if 
$$
\ch^\f(\rho)(\x)\ \geq \ 0 \qquad {\rm for\ all\ } \f-{\rm planes\  } \x {\rm \ tangential \ to\ }  \bo,
\eqno{(\YY.1)}
$$
i.e., for all $\x\in G(\f)$ with span$(\x)\subset T(\bo)$. The boundary $\f$ is {\bf strictly $\f$-convex} if the inequality in (\YY.1) is  strict everywhere on $\bo$.  These conditions are independent of the
choice of defining function $\rho$.

\Theorem{\CC.5} {\sl Let $\O\subset\subset X$ be a compact domain with strictly $\f$-convex boundary, and let $\delta= -\rho$ where $\rho$ is an arbitrary defining function for $\bo$.
Then $-\log \delta:\O\to \bbr$ is strictly \fp outside a compact subset.  In particular,
the domain $\O$ is strictly $\f$-convex at infinity.
}
\medskip

Elementary examples show that the converse of this theorem does not hold in general.
However, there is  a weak partial converse.
\Prop{\CC.8}  {\sl Let $\O\subset\subset X$ be a compact domain with smooth boundary.
Suppose $\f$ is parallel and consider the function $\d = {\rm dist}(\bullet, \bo)$.
If $-\log \d$ is strictly \fp near $\bo$, then $\bo$ is $\f$-convex.}

\medskip

We note that boundary convexity can be  interpreted geometrically as follows.
Let $II$ denote the second fundamental form of the hypersurface $\bo$ oriented by the 
outward-pointing normal.  Then $\bo$ is \fc if and only if trace$(II\bigr|_{\x})\leq 0$
for all $\f$-planes $\x$ which are tangent to $\bo$.

\medskip
\centerline
{$\f$\AAA -P\BBB OSITIVE \AAA C\BBB URRENTS.}\smallskip
\medskip
Our \fp functions are intimately related to the study of $\f$-positive currents introduced
in [HL$_3$]. We recall that a $p$-dimensional current $T$ is called {\bf $\f$-positive }
if it is representable by integration and its generalized tangent $p$-vector
$$
\overrightarrow{T} \in {\rm Convex Hull }(G(\f))\qquad \|T\|-a.e.
$$
where $\|T\|$ denotes the total variation measure of $T$. Examples include $\f$-submanifolds
and, more generally, rectifiable $\f$-currents. By Almgren's   Theorem [A] we know that 
rectifiable $\f$-currents $T$ with $dT=0$ are regular, that is, given by integration over 
$\f$-submanifolds with positive integer multiplicities, outside a closed subset of Hausdorff
dimension $p-2$.

$\f$-Positive currents generalize the positive currents in complex geometry, and $d$-closed
recitifiable $\f$-currents generalize positive holomorphic chains.

Recall from above that if $T$ is a $\f$-positive current (with compact support), then
$$
\supp\, T\ \subset \ \wh{\supp ( d T)} \cup \Core(X).
\eqno{(\YY.2)}
$$
In particular, if $dT=0$, then $\supp T \subset   \Core(X)$, and if $X$ is strictly $\f$-convex ($\Core(X)=\emptyset$), then there exist no $d$-closed $\f$-positive currents with compact support on $X$.

In Section \DD \ we review the known facts concerning  $\f$-positive currents.
These include compactness theorems, regularity theorems, mass-minimizing properties
and dual characterizations.

\medskip
\centerline
{\AAA S\BBB UPERHARMONIC \AAA $	\f$-C\BBB  URRENTS.}\smallskip
\medskip

Assume for now that the calibration  $\f$ is parallel, and consider the adjoint of the operator $dd^\f:\ce^0(X)\to \cd^p(X)$ which can be written as 
$$
\partial_\f\partial:\cd'_p(X) \ \arr\ \cd'_0(X) 
$$
where $\partial:\cd'_p(X) \to \cd'_{p-1}(X) $ denotes the usual adjoint of $d:\ce^{p-1}(X)\to \ce^p(X)$ and 
$\partial_\f:\cd'_{p-1}(X) \ \arr\ \cd'_0(X) $ is the adjoint of $d^\f$, 
defined by $(\partial_\f R)(f)\equiv R( d^\f f)$.

Positive currents $T$ with the property: $\partial_\f\partial T\leq 0$, i.e.,  $\partial_\f\partial T$ is a non-positive measure, satisfy a version of (\YY.2) above.

\Lemma{\EE.2} {\sl Suppose $T$ is a $\f$-positive current with compact support on $X$ which satisfies 
$\partial_\f\partial T\leq 0$ outside a compact subset $K\subset X$. Then
$$
\supp\, T\ \subset \ \wh K \cup \Core(X).
$$
In particular, if $\partial_\f\partial T\leq 0$ on $X$, then  $\supp \,T \subset \Core(X)$.}\medskip

Another consequence is the following.  Suppose $X$ is strictly  $\f$-convex. If
$T$ is   $\f$-positive current with 
 $ \partial_\f\partial T \leq 0$ on $X-K$, then $\supp\, T\subset \wh K$.
In fact, it turns out that  the  points   $ x\in \wh K$ can be characterized in terms
of certain $\f$-positive currents $T$ which satisfy $\partial_\f\partial T = -[x]$ in  $X-K$.

\medskip
\centerline
{\AAA D\BBB UVAL - \AAA S\BBB IBONY \AAA D\BBB UALITY.}\smallskip
\medskip

Points in the $\f$-convex hull of a compact set $K\subset X$ have a useful characterization
in terms of $\f$-positive currents and certain Poisson-Jensen measures.  The following
 results generalize work of Duval and Sibony in the complex case. 
 They remain valid (as does Lemma \EE.2 above) when $\f$ is not parallel if  the operator $\partial_\f\partial$ is replaced with $\ch_\f$.

Let $K\subset X$ be a compact subset and $x$ a point in $X-K$.  By a {\bf Green's current
for } $(K,x)$ we mean a $\f$-positive current $T$ which satisfies
$$
\partial_\f\partial T = \mu-[x]
\eqno{(\YY.3)}
$$
where $\mu$ is a probability measure with support on $K$ and $[x]$ denotes the Dirac measure
at $x$. In this case $\mu$ is called a {\bf Poisson-Jensen} measure for $(K,x)$.
By the remarks above we see that $x\in \wh K$.  In fact we have the following.

\Theorem{\EE.8}  {\sl  Suppose $\f$ is parallel and $X$ is strictly $\f$-convex.
Let $K\subset X$ be a compact subset and $x\in X-K$.  Then there exists a Green's current
for $(K,x)$ if and only if $x\in \wh K$.}

\medskip
We note that if $M\subset X$ is a compact $\f$-submanifold with boundary, and if $G_x$ 
is the Greens function for the riemannian laplacian on $M$ with singularity at $x\in M-\partial M$,
then $\partial_\f\partial(G_x [M]) = \mu-[x]$ for a probability measure $\mu$ on $\partial M$.

\medskip
As a application we obtain the following approximation result.
A domain $\O\subset X$ is said to be {\bf \fc relative to $X$} if 
$K\subset\subset \O \ \Rightarrow  \wh{K}_X\subset\subset \O$.
\Prop{\EE.16} {\sl  Suppose $\f$ is parallel and $X$ is strictly $\f$-convex.
An open subset $\O\subset X$ is  \fc relative to $X$ if and only if 
$\PSH(X,\f)$ is dense in $\PSH(\O,\f)$.}

\medskip

\bigskip
\centerline
{\AAA $\f$-F\BBB REE \AAA S\BBB UBMANIFOLDS  \BBB AND \AAA S\BBB RICTLY
$\f$-\AAA C\BBB ONVEX \AAA S\BBB UBDOMAINS.}\smallskip
\medskip

We next examine the analogues in calibrated geometry of the totally real
submanifolds in complex analysis. Using them we show how to construct strictly \fc manifolds in enormous families with every
topological type allowed by Morse theory.

Let $(X,\f)$ be any fixed calibrated manifold.
A closed submanifold $M\subset X$ is called {\bf $\f$-free} if there are no $\f$-planes
tangential to $M$, i.e.,    no $\x\in G(\f)$ with $\span \x \subset TM$.

Note that $M$ is automatically $\f$-free if it is {\bf $\f$-isotropic}, that is, if $\f\bigr|_M\equiv 0$. . 

Any submanifold of dimension $< p$ is $\f$-free.  Furthermore, generic local submanifolds
of dimension $p$ are $\f$-free.  In some geometries this also holds for certain dimensions $>p$.

In a K\"ahler geometry  the $\o$-free submanifolds are exactly those which are
{\sl totally real} (the tangent spaces contain no complex lines).
In Special Lagrangian geometry,  the $\f$-free submanifolds include the  complex submanifolds
of all dimensions.

\Theorem {\FF.2}  {\sl Suppose $M$ is a closed submanifold of $(X,\f)$ and let 
$\dis(x) \equiv  {\rm dist}(x,M)^2$ denote   the square of the distance to $M$.  Then
$M$ is \totr \ if and only if 
the function
 $\dis$ is strictly \fp at each point in $M$ (and hence in a neighborhood of $M$).
}\medskip

The existence of \totr \  submanifolds insures the existence of many strictly \fc domains in $(X,\f)$.

\Theorem{\FF.4}  {\sl  Suppose $M$ is a \totr\ submanifold of $(X,\f)$.
  Then there exists a fundamental neighborhood system $\cf(M)$ of $M$ such that:
  
  \smallskip
  
\quad  (a)\  $M$ is a deformation retract of each $U\in \cf(M)$.

  \smallskip
  
\quad  (b)\  Each neighborhood $U\in \cf(M)$ is strictly \fc.

  \smallskip
  
\quad  (c)\  \ ${\cal PSH}(V,\f)$ is dense in  ${\cal PSH}(U,\f)$ if $U\subset V$ and $V,U\in\cf(M)$.

  \smallskip
  
\quad  (d)\  Each  compact set $K\subset M$ is ${\cal PSH}(U,\f)$-convex for each $U\in \cf(M)$.
  }

\medskip

This result provides rich families of strictly convex domains.  
Note that neighborhoods $U\in \cf(M)$ include the sets  
$\{x: \dist(x, M)<\e(x)\}$ for positive functions $\e\in C^\infty(M)$ which die arbitrarily rapidly at infinity.
As noted,   any submanifold of  dimension $<p$ is $\f$-free. Furthermore, 
any submanifold of a $\f$-free submanifold is again $\f$-free. 

For example if $X$ is a Calabi-Yau manifold with Special Lagrangian calibration $\f$,
then any complex submanifold $Y\subset X$ is $\f$-free, as is any smooth submanifold
$A\subset Y$. The topological type of such manifolds $A$ can be quite complicated.

This construction can be refined even further by replacing $A\subset Y$ with an arbitrary closed
subset.
It turns out that  the following two classes of subsets of $(X,\f)$:
\smallskip

\qquad (1)\ Closed subsets $A$ of \totr submanifolds.
\smallskip

\qquad (2)\ Zero sets  of non-negative strictly \fp\ functions $f$.
\smallskip

\noindent
are essentially  the same. This is given in Propositions \FF.7 and \FF.8
One also has:
\Prop{\FF.12} {\sl Let $f$ be a non-negative, real analytic function on $(X,\f)$ and consider the real analytic subvariety $Z\equiv \{f=0\}$.  If $f$ is strictly \fp at points of $Z$, then each stratum of $Z$
is \totr.}

The results above generalize   work of Harvey-Wells [HW$_{1,2}$] in the complex case.

 \bigskip

\centerline
{\AAA B\BBB OUNDARIES \BBB OF \AAA $\f$-S\BBB UBMANIFOLDS.}\smallskip
\medskip

A very natural question in calibrated geometry is the following:  Given a compact
oriented submanifold $\G\subset X$ of dimension $p-1$, when does there exist
a $\f$-submanifold $M$ with $ \partial M=\G$? A companion question is: Given 
a compactly supported current $S$ of dimension $p-1$ in $X$, when does there exist
a $\f$-positive current $T$ with $\partial T=\G$? For this second question there is a complete answer
when   $X$ is strictly \fc and $\f$ is exact.

\Theorem {\HH.1} {\sl   Consider  a current
$S\in \ce'_{p-1}(X)$. Then $S=\partial T$  for some  $\f$-positive current $T\in \ce'_p(X)$
if and only if }
$$
\int_S \, \a\ \geq \ 0\qquad\ \ {\rm for\ all\ } \a\in \ce^{p-1}(X) \ \ \ {\rm such\ that\ } d\a \ {\rm is\ }
\lpp{\rm -positive}
$$
\Note{} $\lpp$-positive means that $d\a(\x) \geq 0$ for all $\x\in G(\f)$.

There is a similar result for compact calibrated manifolds $(X,\f)$ with
no condition on $\f$.  The result uses $\f$-quasiplurisubharmonic functions --
those which satisfy the condition that $dd^\f f+\f$  is $\lpp$-positive.  See [HL$_4$] for the K\"ahler case.
\medskip

\vfill\eject
\centerline
{\AAA $\f$-F\BBB LAT  \AAA H\BBB YPERSURFACES.}\smallskip
\medskip

In Section \II \ we expand the notion of $\f$-pluriharmonic functions to include functions $f$
which are $\f$-pluriharmonic modulo the ideal generated by $df$.
In most interesting geometries these functions are characterized by the fact
that their level sets are {\bf $\f$-flat}, i.e., the trace of the second fundamental form
on all tangential $\f$-planes is zero.  These functions are important for the boundary problem.
If $f$ is such a function defined in a neighborhood of a compact $\f$-submanifold with boundary  
 $M\subset X$, then
$$
\inf_{\partial M} \ \leq f(x) \ \leq \ \sup_{\partial M} f \qquad\ \ {\rm for \ \ } x\in M.
$$

\bigskip
\centerline
{\AAA G\BBB ENERALIZED \AAA H\BBB ODGE  \AAA M\BBB ANIFOLDS.}\smallskip
\medskip

 In Section \GG\ we discuss analogues of  Hodge manifolds in the general calibrated setting.
We also examine various analogues of the Hodge Conjecture in these spaces.

\medskip

The paper is organized  in the order presented above.  Several of the sections
have appendices which contain examples or discussions of side issues.
They can be skipped when reading the paper the first time.

We mention that the operator $d^\f$ has been independently found by M. Verbitsky 
[V] who studied the generalized K\"ahler theory (in the sense of Chern) on G$_2$-manifolds.
The authors would like to thank Robert Bryant for useful comments and conversations
related to this paper.

\vfill\eject


\centerline{\bf \AA. Plurisubharmonic Functions}
\medskip

Suppose $\phi$ is a calibration on a riemannian manifold $X$.  
The  $\phi$-Grassmann bundle, denoted $G(\phi)$, consists of the unit simple
vectors $\xi$ with $\phi(\xi)=1$, i.e., the $\phi$-planes.  An oriented submanifold $M$ is  a
{\sl $\phi$-submanifold}, or is {\sl calibrated by $\phi$}, if the
oriented unit tangent space ${\oa T}_x M$ lies in $G_x(\phi)$ for each
$x\in M$, or equivalently, if $\phi$ restricts to $M$ to be the volume form
on $M$.  Let $n = \dim X$ and $p = {\rm degree}(\phi)$. 

\Def {\AA.1} The {\bf $\df$-operator} is defined by $$\df f \ =\ \nabla f
\hk      \phi$$ for all smooth functions $f$ on $X$. 

\medskip
\noindent
Hence
$$
\df : \ce^0(X)\ \arr\ \ce^{p-1}(X) \and d \df :\ce^0(X)\ \arr\ \ce^{p}(X)
$$
where   $\ce^p(X)$ denotes the space of $C^\infty$ $p$-forms on $X$. 
This $d d^\f$ operator  provides a way of defining
plurisubharmonic functions in   calibrated geometry when the calibration $\f$ is parallel.

If  $\o$ is a K\"ahler form on a complex manifold, then
$
d^\o \ =\ d^c\ =\ -J\circ d
$
is the conjugate differential. Thus, the $dd^\f$-operator generalizes the 
$dd^c$-operator in complex geometry.

\Def {\AA.2} Suppose $\nabla \phi=0$. A  function $f\in C^\infty(X)$
is {\bf $\phi$-plurisubharmonic} if 
$$
(d \df f)(\xi)\ \geq 0 \qquad {\rm for\  all } \ \ \xi \in G(\phi).
$$
The set of such functions will be denoted $\cp\cs\ch(X,\phi)$.
If $(d \df f)(\xi)>0$ for all $\xi \in G(\phi)$, then $f$ is {\bf strictly  $\phi$-plurisubharmonic}. 
If $(d \df f)(\xi)=0$ for all $\xi \in G(\phi)$, then $f$ is  {\bf $\f$-pluriharmonic}.

\Remark {\AA.3}  If $\f$ is not parallel, we define \fp functions by replacing $dd^\f f$,
in Definition \AA.2, with  $dd^\f f - \nabla_{\nabla f} \f$.  This modified $dd^\f$-operator is discussed in detail in Section \ZZ. Note that the difference $\nabla_{\nabla f} \f$ is a first order operator.

\medskip

The next result justifies the use of the word plurisubharmonic in the
context of  a $\phi$-geometry.  A calibration $\f$  is {\bf integrable} if for
each point $x\in X$ and each $\x\in G_x(\f)$ there exists a $\f$-submanifold $M$ through $x$ with 
$\oa{T_xM} =\x$.

\Theorem{\AA.4} {\sl  Let $(X,\f)$ be any calibrated manifold.  If a function  $f\in C^\infty(X)$ is $\phi$-plurisubharmonic, then the restriction of $f$ to any $\phi$-submanifold $M\subset X$ is subharmonic in the induced  metric.  If $\f$ is integrable, then the converse holds.}
\medskip

Theorem \AA.4  is an immediate consequence of the formula
$$
\left(    dd^\f f - \nabla_{\nabla f} \f   \right)\bigr|_M\ =\ (\D_M f )\, {\rm vol}_M
\eqno{(\AA.1)}
$$
This formula follows  from the three equations (\ZZ.7), (\ZZ.12), and (\ZZ.15), proved
below, and the fact that $\f$-submanifolds are minimal submanifolds.  We continue for the moment to present results whose proofs will be given  in  Section \ZZ.

The $\phi$-plurisubharmonic functions enjoy many of the
useful properties of their classical cousins in complex analysis.
Here are some basic examples.

\Prop {\AA.5} {\sl Let $f,g\in C^\infty(X)$ be
$\phi$-plurisubharmonic.\smallskip

(i)\ \  If $\psi\in C^\infty(\bbr)$ is convex and increasing, then $\psi \circ
f$ is $\phi$-plurisubharmonic.\smallskip

(ii) \ The function $\log(e^f+e^g)$ is $\phi$-plurisubharmonic.
\smallskip

(iii) The decreasing sequence of functions $h_n\equiv {1\over n}\log \left(e^{nf}+e^{ng}\right)$ of 
\fp functions approximates $\max\{f,g\}$.  More precisely, }
$
h_n -{1\over n}\log 2 \ \leq \ \max\{f,g\} \ \leq \ h_n.
$

\medskip
Another important elementary property is given in the next proposition.  A $p$-form 
$\f$ is said to {\bf  involve all the variables} if $\zeta\hk \f\neq 0$ for all non-zero tangent vectors $\z$. 

\Prop{\AA.6}  {\sl  The $dd^\f$-operator is (overdetermined) elliptic if and only if the calibration involves all the variables.  Its  symbol is }
$$
\s_{\z}(d d^\f)\ =\ \z\wedge(\z\hk\f) \qquad \ \  {\rm for \ }\    \z\in T^*_xX \cong T_xX.
\eqno{(\AA.2)}
$$
\pf  The computation of (\AA.2) is straightforward.
By definition the operator  $dd^\f$ is (overdetermined) elliptic if  $\s_{\z}(d d^\f)$ is injective
 (i.e., $\ne 0$)  for all
$\z\neq 0$.   Observe  that 
$\z\wedge(\z\hk\f)=0$ if and only if $\z\hk \f=0$ since 
$\langle \z\wedge(\z\hk\f), \f\rangle = | \zeta\hk \f |^2$.  Hence,  
$ \z\wedge(\z\hk\f) \ne0$ for all $\z\ne0$ if and only if
$\z\hk \f\ne0$ for all $\z\ne0$. \qed

\Def{\AA.7}  The {\bf $\phi$-Laplacian}, $\Delta_{\phi}$,  on a function $f\in C^\infty(X)$ is
defined by
$$
\Delta_{\phi} f \ \equiv\ \langle d \df f, \phi\rangle
$$
 The symbol of $\D_\f$ is the quadratic form  $|\zeta \hk \f|^2$ for  $\zeta\in T^*_xX
\cong T_xX$.      In particular,   the $\phi$-Laplacian is (determined) elliptic  if and only if the calibration involves all the variables.
\medskip

Definition \AA.7 is easily extended to include arbitrary distributions $f$  on $X$ by requiring 
$(dd^\f f, \a\otimes *1)\geq 0$ for all smooth, compactly supported sections $\a$ of the bundle $\Lambda_pTX$
which are positive linear  combinations of elements in $G_x(\f)$ at each point $x$.

\Prop{\AA.8} {\sl    Suppose $\f$ is a calibration  which can be written  as a positive linear combination
of $\f$-planes at each point $x\in X$.  If a distribution $f$ is $\f$-plurisubharmonic, then $f$ is $\D_\f$-subharmonic,
i.e.,  $\D_\f f$ is a non-negative measure.}
\medskip

\Remark{}    If $\f$ is a calibration which satisfies the hypotheses of both Propositions \AA.6 and \AA.8, then all the  classical results 
concerning subharmonic functions with respect to the scalar elliptic operator $\D_\f$ can be brought to bear on \fp distributions.  For example, each  \fp distribution is in $L^1_{\rm loc}$ 
(locally Lebesgue integrable) and represented by an upper-semicontinuous function taking
values in   $[-\infty, \infty)$.   However, because of the geometric emphasis in this paper and the
desire to keep technical considerations to a minimum, only $C^\infty$ \fp functions will be considered.
Exceptions will occur as remarks.
\medskip

Suppose $\f$ is a calibration for which $\L(\f) \equiv \span G(\f) \subset \Lambda^pT^*X$
is a subbundle.   Then {\sl for any  \fp distribution $f$,   the $\L(\f)$-component   of $dd^\f f$  has measure-coefficients}. This  is proved in Proposition \DD.19  below.

\medskip

\Remark{{\bf (The Abundance of \fp Functions)}}   We shall see in the next section that
any convex function on the riemannian manifold $X$ is automatically \fp.  However, in the
euclidean case $(\bbr^n, \f)$ with $\f$ parallel,  there are many   strictly \fp quadratic functions which are not convex.  (See Remark 2.9.)  It follows that in small neighborhoods of a point
on any calibrated manifold,  such functions exist.  Of course the strictly \fp functions form an open cone
in $C^\infty(X)$ for any calibrated manifold $(X,\f)$.



\vskip .3in

\centerline{\bf  Appendix:   Pluriharmonic Functions} 

\medskip
 While  \fp functions 
are abundant, the $\f$-pluriharmonic functions are often quite scarce.  To illustrate this 
phenomenon we shall sketch some of the basic facts in the ``classical'' cases.

To begin we note that for some calibrations $\f$, one has that:
$$
dd^\f f\ =\ 0 \qquad {\rm if\ and\ only\ if\ }\qquad (dd^\f f)(\x) \ =\ 0 \quad {\rm for\ all\  } \x\in G(\f)
\eqno{(\AA.3)}
$$
while for others this is not true. It is the right hand side that defines pluriharmonicity.
If (\AA.3) holds and the map $\BM_\f$ is everywhere injective (as in Example \AA.11), then the only pluriharmonic functions are the affine functions, i.e., the functions with parallel gradient.
Note that if $f$ is affine, then $\nabla f$ splits the manifold locally as a riemannian 
product $X=\bbr\times X_0$.  

\Ex{\AA.8. {\bf (Complex geometry)}}  Let $\o$ be a K\"ahler form on a complex manifold
$X$.  Then
$
d^\o = d^c
$
is the conjugate differential, $dd^cf$ is the complex hermitian Hessian of $f$, $G(\o)$ is the grassmannian of complex lines, and the statement (\AA.3) is valid. In particular, the 
$\o$-pluriharmonic functions are just the classical pluriharmonic functions on $X$.

\Ex{\AA.9. {\bf (Special Lagrangian geometry)}}   Consider the special Lagrangian calibration
$\f={\rm Re}( dz)$ on $\bbc^n$.  Let $Z_{ij}$ denote the bidegree $(n-1,1)$ form obtained
from $dz = dz_1\wedge\cdots\wedge dz_n$ by replacing $dz_i$ with $d{\bar z}_j$  (in the $i$th position).    
A short calculation shows that
$$
dd^\f f \ =\ 2 {\rm Re} \biggl\{   \sum_{i,j=1}^n{\partial^2 f \over \partial {\bar z}_i \partial {\bar z}_j}
Z_{ij}\biggr\}
+(\D f) {\rm Re}( dz)
\eqno{(\AA.4)}
$$
For this calibration one can show that (\AA.3) is valid.  Consequently, Lei Fu [Fu] has described 
all the $\f$-pluriharmonic functions.

\Prop{\AA.10} {\sl Let $f$ be a special Lagrangian pluriharmonic function defined locally on $\bbc^n$, $n\geq 3$.  Then $f=A+Q$ where $A$ is affine and $Q$ is 
a traceless hermitian quadratic function.}

\pf  If $dd^\f f=0$ and $n\geq 3$ (so that $Z_{ij}$ and ${\overline Z}_{ij}$ are of different bi-degrees), then(\AA.4) implies that 
${\partial^2 f \over \partial {\bar z}_i \partial {\bar z}_j} = 0$ for all $i,j$.  Therefore, all third partial derivatives if $f$ are zero.  For polynomials of degree $\leq 2$ the result is transparent from (\AA.4).\qed

\vfill\eject

\Ex{\AA.11. {\bf (Associative, Coassociative and Cayley geometry)}}  
Consider one of the calibrations:\smallskip

1. (Associative) \qquad $\f(x\wedge y\wedge z) \ =\ \langle x, yz\rangle$ \qquad for $x,y,z \in {\rm Im} \bbo$

2. (Coassociative) \qquad\ $\psi \ =\ *\phi$  \qquad \qquad\qquad\qquad \ \ \ on ${\rm Im} \bbo$

3. (Cayley) \qquad $\Phi(x\wedge y\wedge z\wedge  w) \ =\ 
\langle x, y\times z\times w\rangle$ \qquad for $x,y,z, w \in  \bbo$

\smallskip\noindent
where $\bbo$ denotes the octonions.  As in the special Lagrangian case
one can show that (\AA.3) is valid for each of these calibrations.   Furthermore, an application
of representation theory shows the maps $\BM_\f, \BM_{\psi}$ and $\BM_\Phi$ are injective.
 These calculations carry over to manifolds with  $G_2$ or   Spin$_7$-holonomy  to establish the following. 

\Prop{\AA.12}  {\sl Let $X$ be a manifold with holonomy contained in  $G_2$ or Spin$_7$
and  having dimension 7 or 8 respectively.  Suppose $\f$ is a parallel calibration
on $X$  of one of the three types above.  Then every  $\f$-pluriharmonic function
on $X$ is affine.  Moreover,  if the holomony is exactly $G_2$ or Spin$_7$,
 every $\f$-pluriharmonic function is constant.  }
 
 \pf
 The first assertion follows because (\AA.3) is valid and the $\BM$-maps are injective.
 The second follows because any non-constant affine function on $X$ would reduce
 its  holonomy to a subgroup of $1\times SO_{n-1}$. \qed

\Ex{\AA.13. {\bf (Quaternionic-K\"ahler  geometry)}}  
Let $\bbh$ denote the quaternions and consider $\bbh^n$ as a right-$\bbh$
vector space.  Each of the complex structures $I, J, K$ (right multiplication by $i,j,k$) determines a
 K\"ahler form $\o_I, \o_J, \o_K$ respectively.  The 4-form
 $$
 \Psi\ \equiv\  {1\over6}(\o_I^2 + \o_J^2 + \o_K^2)
 \eqno{(\AA.5)}
 $$
on $\bbh^n\equiv \bbr^{4n}$ is a calibration with $G( \Psi)$ consisting of the oriented quaternion
lines in $\bbh^n$.  In this case, $d d^ \Psi f \equiv 0$ if and only if $\Hess f \equiv 0$.  However,
the assertion (\AA.3) is not valid in this case, and in fact there is a rich family of $ \Psi$-pluriharmonic functions.  For example, if $f$ is $\o_I$-pluriharmonic,  then $f$ is $ \Psi$-pluriharmonic.
Hence, so is any $\o$-pluriharmonic $f$ where $\o=a\o_I+b\o_J+c\o_K$ with $a^2+b^2+c^2=1$.

It is well known that the only $ \Psi$-submanifolds in $\bbh^n$ are the
affine quaternion lines.

Of course the calibration (\AA.5) exists on any quaternionic K\"ahler manifold,
i.e.,  one with ${\rm Sp}_n \cdot {\rm Sp}_1$-holonomy.  (See [GL] for examples.)   With this 
full  holonomy group it seems unlikely that there are many
$\Psi$-pluriharmonic functions. However, if the holonomy is contained in SP$_n$, they
exist in abundance as seen in the next example.

\Ex{\AA.14. {\bf (Hyper-K\"ahler  manifolds)}}  Let $(X, \o_I,\o_J,\o_K)$ be a hyper-K\"ahler manifold.
Then $X$ carries  several parallel calibrations. There are, of course, the K\"ahler forms
$\o=a\o_I+b\o_J+c\o_K$ with $a^2+b^2+c^2=1$, and two others of particular interest.

(1) \ Let $\Psi = {1\over6}(\o_I^2 + \o_J^2 + \o_K^2)$.  Then as in Example \AA.13, {\sl any 
$\o$-pluriharmonic function is $\Psi$-pluriharmonic}.  Hence, the sheaf of $\Psi$-pluriharmonic
functions is quite rich. On the other hand there are precious few $\Psi$-submanifolds.

(2)\  Consider the generalized Cayley form 
$
\Xi \ \equiv\ {1\over2}(\o_I^2-\o_J^2-\o_K^2).
$
For this calibration there exist no interesting pluriharmonic functions, but there
are many $\Xi$-submanifolds (cf. [BH]).

\Ex{\AA.15. {\bf (Lie  group geometry)}}  Let $G$ be a compact simple Lie group with Lie algebra
$\gerG$, defined as the set of left-invariant vector fields on $G$.  Consider the fundamental 
3-form $\f$ on $G$ defined by $\langle x, [y,z]\rangle$ and normalized to have comass one.
Calculations indicated that in all but a finite number of cases non-constant pluriharmonic functions
do not exists, however there are $\f$-submanifolds, namely the ``minimal'' SU$_2$-subgroups
(cf. [B]).

\Ex{\AA.16. {\bf (Double point  geometry)}}  Let $\f = dx_1\wedge\cdots\wedge dx_n +
 dy_1\wedge\cdots\wedge dy_n$  in $\bbr^{2n}$ for $n\geq 3$.  The only $\f$-planes are those parallel to the $x$ or $y$ axes.  An easy calculation shows that $dd^\f f=0$ 
 if and only if $f(x,y)=g(x)+h(y)$ for harmonic functions $g$ and $h$. However,  a function 
 $f(x,y)$ is $\f$-pluriharmonic if and only if it is harmonic in $x$ and $y$ separately.

\vfill\eject



\centerline{\bf 2. \ The $\f$-Hessian.}\medskip

In this section we prove the assertions made in \S 1.  The arguments will involve ideas and notation
important for the rest of the paper.  We will end the section with a generalization of Theorem 1.4  to submanifolds which are {\sl $\f$-critical}.

Recall that the {\bf Hessian} (or second covariant derivative)  of a smooth function $f$ on a riemannian
manifold $X$ is defined on tangent vector fields
$V,W$ by
$$
\Hess(f)(V,W) \ \equiv\ V( W  f) - (\nabla_V W) f
\eqno{(\ZZ.1)}
$$
where $\nabla$ denotes the riemannian connection.
Note that $V( W  f) - (\nabla_V  W) f = V(\langle  W, \nabla f\rangle) - \langle \nabla_V  W, \nabla f\rangle = \langle W, \nabla_V(\nabla f)\rangle$ so that at a point $x\in X$, the Hessian is the symmetric 
2-tensor, or the symmetric linear map of $T_xX$ given by
$$
\Hess(f)(V) \ =\   \nabla_V(\nabla f) 
\eqno{(\ZZ.2)}
$$

Let $V$ be a real inner product space 
Given an element $\f\in \L^p V^*$, we define a linear map, central to this paper, 
$$
{\BM}_{\f} : \End(V)\ \arr\  \L^p V^*
\eqno{(\ZZ.3)}$$
by
$$
{\BM}_{\f} (A) \ \equiv \ {\Der}_{A^t}(\f)
$$
where ${\Der}_{A^t}$ denotes the extension of the transpose $A^t:V^*\to V^*$ 
to ${\Der}_{A^t}:\L^pV^*\to \L^pV^*$ as a derivation.

\Note{} Recall that the natural inner product on $\End(V)$ is given by:
$$
\langle A,B   \rangle \ =\ \tr A B^*  \qquad {\rm for }\  A,B \in \End(V)
$$
Using this inner product we have  the adjoint map
$$
{\BM}_{\f}^* : \L^p V^*\ \arr\   \End(V)
\eqno{(\ZZ.4)}
$$
which will also be important.

\Def{\ZZ.1}  The {\bf $\f$-Hessian} of a function $f\in C^\infty(X)$ is the $p$-form $\CH^\f(f)$
defined by letting  the symmetric endomorphism    $\Hess f$ act on $\f$ as a derivation, i.e., 
$$
\CH^\f(f)   \ \equiv \    \Der_{\Hess^t f}( \f).
\eqno{(\ZZ.5)}
$$
In terms of the bundle map ${\BM}_{\f}  : \End(TX)\to \L^pT^*X$,
$$
\CH^\f(f)   \ \equiv \  {\BM}_{\f} (\Hess f)
\eqno{(\ZZ.6)}
$$
is the image of the Hessian of $f$.
\medskip

The second order differential operators $dd^\f$ and $\CH^\f$ differ by a pure 
first order operator. This is the first of the three equations needed to prove
Theorem \AA.4.

\Theorem{\ZZ.2}  {\sl  If $\f$ is a closed form on $X$, then}
$$
\CH^\f(f) \ =\ dd^\f f- \nabla_{\nabla f} \f
\eqno{(\ZZ.7)}
$$
\pf
By (\ZZ.2)  we have  $(\Hess f)(V) = \nabla_V \nabla f = [V, \nabla f] +\nabla_{\nabla f} V$, i.e.,
$$
\Hess f\ =\   -\cl_{\nabla f} + \nabla_{\nabla f}
$$
as operators on vector fields ($\cl$ = the Lie derivative).
The right  hand side of this formula has a standard extension to all tensor fields
as a derivation that commutes with contractions.  It is zero on functions, that is, it is
a bundle endomorphism whose value on $T^*X$ is minus the transpose of its value
on $TX$.  In particular, we find that
  ${\Der}_{\Hess^t f} = \cl_{\nabla f} - \nabla_{\nabla f}$ on $p$-forms, i.e., 
$$
\CH^\f(f) \ =\ \cl_{\nabla f}(\f) - \nabla_{\nabla f}\f
\eqno{(\ZZ.8)}
$$
Finally, since $d\f=0$, the classical formula relating $\cl,d$ and $\hk$ gives
$$
\qquad\qquad\qquad\qquad\qquad\qquad
dd^\f f \ =\ d(\nabla f\hk \f)\ =\ \cl_{\nabla f}(\f)  \qquad\qquad\qquad\qquad\qquad
\vrule width5pt height5pt depth0pt
$$
\medskip

Many of the nice results for the $dd^\f$-operator continue to hold in the 
non-parallel case after replacing it with the $\f$-Hessian. Perhaps even more importantly,
many properties of the $dd^\f$-operator  in the parallel case can best be understood by
considering the $\f$-Hessian. 

The second formula needed for the proof of Theorem \AA.4 is algebraic in nature,
involving  the bundle map ${\BM}_{\f}  : \End(TX)\to \L^pT^*X$.  Consequently, as before,
we replace $T_xX$ by a general inner product space $V$. If $\x$ is a $p$-plane
in $V$ (not necessarily oriented),  let $P_{\x}:V \to \x$ denote orthogonal projection.
The following, along with its reinterpretations (\ZZ.9$)'$ and (\ZZ.12), is a central result of this paper.

\Theorem{\ZZ.3}  {\sl  Suppose $\f$ has comass one.  For each $A\in \End(V)$,
$$
({\BM}_{\f}  A)(\x) \ =\ \langle  A, P_\x \rangle \qquad {\rm if\ } \x \in G(\f).
\eqno{(\ZZ.9)}
$$
Equivalently,}
$$
({\BM}_{\f}^*)(\x) \ =\ P_\x  \qquad {\rm if\ } \x \in G(\f).
\eqno{(\ZZ.10)}
$$
\medskip
Note that if $e_1,...,e_p$ is an orthonormal basis for the $p$-plane $\x$, then
$$
\langle  A, P_\x \rangle \ =\ \sum_{j=1}^p  \langle  e_j, Ae_j \rangle 
$$
Consequently, it is natural to refer to  $\langle  A, P_\x \rangle $  as the {\bf $\x$-trace of $A$}
and to use the notation 
$$
\tr_\x A\ \equiv\  \langle  A, P_\x \rangle.
$$
In particular, for each $A\in \End(V)$,
$$
({\BM}_{\f}  A)(\x) \ =\ \tr_\x A  \qquad {\rm if\ } \x \in G(\f).
\eqno{(\ZZ.9)'}
$$

Suppose $\x \in G(p,V)\subset \L_p V$ is a unit simple vector.  If
$a,b$ are unit vectors in $V$ with $a\in \span \x$ and $b\perp
\span\x$, then  
$$
b\wedge(a\hk\x)
$$
is called a {\bf  first cousin} of $\x$.  The first cousins of $\x$ form a
basis for the tangent space to the Grassmannian 
$G(p,V)\subset \L_p V$ at the point $\x$. 
Since $\f$ restricted to $G(p,V)$ is a maximum on $G(\f)$, this fact implies the
following  result, which we shall use frequently.

\Lemma{\ZZ.4. {\bf (The First Cousin Principle)}} {\sl  If $\phi\in    \L^p V^*$ has comass one and
$\x \in G(\f)$, then 
$$
\phi(\eta)\ =\ 0
$$
for all first cousins $\eta = b\wedge(a\hk\x)$ of $\xi$.}

\medskip

Note that $\Der_{(b\otimes a)^t} \f =\Der_{a\otimes b} \f = a\wedge(b\hk \f)$ and $\Der_{b\otimes a} \x =  b\wedge(a\hk\x)$
so that if $A=b\otimes a$ is rank one, then
$$
{\BM}_{\f}  (a\otimes b)(\x)\ =\ (\Der_{a\otimes b} \f )(\x)\  = \ \f(\Der_{b\otimes a} \x )
\ =\ \f(b\wedge(a\hk\x))
\eqno{(\ZZ.11)}
$$

\noindent
{\bf Proof of Theorem \ZZ.3.}  Pick an orthonormal basis for $\x$ and extend to an orthornormal basis of $V$.  It suffices to prove (\ZZ.9) when   $A=b\otimes a$
with $a$ and $b$ elements of this basis. It is easy to see that  $\langle  b\otimes a, P_\x \rangle =0$ 
unless $a=b\in\x$, in which case $\langle  a\otimes a, P_\x \rangle =1$.
By equation (\ZZ.11) we have ${\BM}_{\f}  (b\otimes a)(\x) = \f(b\wedge(a\hk\x))$ and
$b\wedge(a\hk\x)=0$ unless $a\in \x$ and either $b\in\x^{\perp}$ or $b=a$.
If $b\in\x^\perp$, then $(b\wedge(a\hk\x)$ is a first cousin of $\x$ and $\f((b\wedge(a\hk\x))=0$ by
the First Cousin Principle.  If $a=b\in\x$, then $b\wedge(a\hk\x)=\x$ and therefore
$\f((b\wedge(a\hk\x))=\f(\x)=1$. \qed

\medskip

Theorem \ZZ.3 has many   consequences.  We mention several. From (\ZZ.9$)' $ we have:

\Cor{\ZZ.5} {\sl  Suppose $(X,\f)$ is a calibrated manifold.  For each function $f\in C^\infty(X)$,}
$$
{\CH}^\f(f)(\x)\ =\ \tr_{\x}(\Hess f)  \qquad {\rm if\ }  \x\in G(\f).
\eqno{(\ZZ.12)}
$$

\medskip
This equation (\ZZ.12) is the second equation needed in the proof of Theorem \AA.4.

\Remark{} Equation  (\ZZ.12) provides an alternative definition of $\f$-plurisubharmonic
(as well as strict $\f$-plurisubharmonic and $\f$-pluriharmonic) functions, which bypasses
the bundle map ${\BM}_{\f} $ and uses only the trace of the Hessian of $f$ on $\f$-planes $\x$.

Another application of Theorem \ZZ.3 is given by:

\Cor{\ZZ.6} {\sl   If $A\in\End(V)$ is skew, then the $p$-form ${\BM}_{\f}  A$ vanishes on $G(\f)$.}
\medskip

Theorem \ZZ.3 can be used to prove Proposition \AA.5.  Note that for $A,B\in \Sym^2(V)
\subset\End(V)$, if   $A\geq 0$, $B \geq 0$, then $\langle A,B\rangle =\tr AB  \geq0$.
Hence for all $\x\in G_p(V)$ one has  $\langle e \otimes e, P_\x\rangle \geq0$, and more generally
 $\langle A, P_\x\rangle \geq0$ whenever $A\geq0$.  Since $df$ and $\nabla f$ are metrically equivalent, 
 $$
 {\BM}_{\f} (\nabla f \otimes \nabla f) = df\wedge ( \nabla f \hk \f) = df \wedge d^\f f.
\eqno{(\ZZ.13)} $$
 Therefore,  Theorem \ZZ.3 has the following consequence.  

\Cor{\ZZ.7} {\sl  For any $f\in C^\infty(X)$, }
$$
(df\w d^\f f)(\x)  \ =\ |\nabla f \hk \x |^2          \  \geq \ 0 \qquad{\sl for\ all\ }  \ \x \in G(\f).
$$
\medskip
\noindent
{\bf Proof of Proposition \AA.5.} We will use Corollary \ZZ.7.  For (i) note  that 
$$
dd^\f(\psi\circ f)\ =\ (\psi'\circ f) \, dd^\f f+ (\psi''\circ f)\, df\w
d^\f f,
$$
which combined with (\ZZ.7) shows that
$$
\ch^\f(\psi\circ f)\ =\ (\psi'\circ f) \,\ch^\f f+ (\psi''\circ f)\, df\w
d^\f f. 
\eqno(\ZZ.14)
$$
For (ii) compute that:
$$
\ch^\f \log(e^f+e^g) \ =\ 
\bigg({{e^f}\over{e^f+e^g}}\bigg)\ch^\f f +
\bigg({{e^g}\over{e^f+e^g}}\bigg)\ch^\f g +
{{e^fe^g}\over{(e^f+e^g)^2}}d(f-g)\w d^\phi(f-g).  
$$
For (iii) set $a=e^f$ and $b=e^g$ in the  inequalities:
$
\max\{a,b\} \ \leq\ \root n \of {a^n+b^n} \ \leq\  2^{1\over n} \max\{a,b\} 
$
and take the log.
\qed
\medskip

Theorem \ZZ.3 can also be used to understand the relationship between convex functions and \fp functions.  A function $f\in C^\infty(X)$ is called {\bf convex} if $\Hess f\geq 0$ at each point,
and it is called {\bf affine} if $\Hess f \equiv 0$ on $X$.  (If $f$ is affine, $\nabla f$ splits $X$ locally as a riemannian product
$\bbr\times X_0$.)

\Cor{\ZZ.8}  {\sl  Every convex function is \fp, and every strictly convex function is strictly \fp  (and every affine function is $\f$-pluriharmonic).}

\medskip

\Remark{\ZZ.9} The converse always fails; there are always \fp functions which are not convex.
To see this, consider first  the euclidean case with $X=V$ and $\f$ parallel.  
Recall that the orthogonal projections $P_e$ onto lines in $V$ generate  the extreme
rays of the convex cone of convex functions (positive semi-definite quadratic forms) in $\Sym^2V\subset \End(V)$.  This cone is self-dual. The projections $P_\x={\BM}_{\f}^*(\x) $ for $\x\in G(\f)$
generate a  proper convex subcone  (in fact a proper convex subcone of the cone generated by orthogonal projections onto $p$-planes).  Hence, by the Bipolar Theorem  there must exist a non-convex quadratic function $Q\in \Sym^2 V$ with $\langle Q, P_\x\rangle \geq 0$ for all $\x \in G(\f)$.  By (\ZZ.9), $Q$ is \fp.\medskip

We now recall some elementary facts about submanifolds.
Given a submanifold $\overline X\subset X$, let $(\bullet)^T$ and 
$(\bullet)^N$ denote orthogonal projection of $T_xX$ onto the tangent and
normal spaces of $\ol X$ respectively. Then the canonical riemannian
connection $\overline \nabla$ of the induced metric on $\overline X$ is
given by ${\ol \nabla }_{ V}  W = (\nabla_V W)^T$ for tangent vector fields
$V,W$ on $\ol X$. The {\bf  second fundamental form} is defined by
$$
B_{V,W} \ \equiv \ (\nabla_V W)^N \ =\ \nabla_V W -\ol\nabla_V W.
$$
This is a symmetric bilinear form on $T\ol X$ with values in the normal
space. Its trace $H={\rm trace}\, B$ is the {\bf  mean curvature vector field}
of $\ol X$, and $\ol X$ is called a {\bf  minimal submanifold } if $H\equiv0$. 
Finally, let $\overline{\Delta}$ denote the Laplace-Beltrami operator   on $\ol X$
and  $\ol{\Hess}$ denote the Hessian operator on  $\ol X$.  The proof of the following 
is straightforward.
 $$
\ol \Hess(f)(V,W) = \Hess(f)(V,W) - B_{V,W} \cdot f
$$ 
for tangent vectors $V,W\in T \ol X$.  Taking the $T\ol X$-trace yields:
$$
\ol{\D} \ =\ \tr_{T\ol X} \Hess f  - H(f).
$$

With a change of notation, this is the   final formula needed to prove Theorem \AA.4.

\Prop{\ZZ.10}  {\sl  Suppose $M$ is a $p$-dimensional submanifold of $X$ with mean curvature vector field
$H$.  Then for each $f\in C^\infty(X)$, }
$$
{\D_M} f \ =\ \tr_{TM} \Hess f  - H(f)  \qquad  {\ \rm on\ } M
\eqno{(\ZZ.15)}
$$

\Cor{\ZZ.11} {\sl Suppose $M$ is a $\f$-submanifold of $X$. Then}
$$
\ch^\f(f)\bigr|_M\ =\ (\D_Mf) {\rm vol}_M
\eqno{(\ZZ.16)}
 $$
\pf Combine (\ZZ.12) and (\ZZ.15) with the fact that $H=0$.\qed\medskip

Combining this with   (\ZZ.7) gives equation (\AA.1) and proves
Theorem \AA.4.
This completes the proof of all  the results in \S \AA \  except Proposition \AA.8.

 \medskip
\noindent
{\bf Proof of Proposition \AA.8.}   
By definition, a distributional section of a vector bundle $E\to X$ is a continuous linear functional on the space of smooth compactly supported sections of
$E^*\otimes \L^nT^*X$, or equivalently, on the space of $\wt s \equiv s\otimes *1$ for 
$s\in \G_{\rm cpt} E^*$.

Suppose $f$ is a \fp distribution, that is, 
$(dd^\f f, \wt\a)\geq 0$
for all smooth compactly supported sections $\a$ of the bundle $\L_pTX$ which are positive linear 
combinations of $\f$-planes at each point.  
Let $\wt\f$  denote the section of $\L_pTX$ corresponding to $\f$ under the metric equivalence   $\L_p\cong \L^p$.  Set $\wt{\a}=g\wt \f$ with $g\geq 0$ a smooth compactly supported function.
By hypothesis $0\leq (dd^\f f, g \wt \f)$, and we claim that $ (dd^\f f, g \wt \f) = (\D_{\f}  f, \wt g)$ where $\wt g \equiv g(*1)$. To see 
this  it suffices to consider the case where $f$ is smooth, where one has $(dd^\f f, g\wt \f) = (\langle dd^\f f, \f\rangle, \wt g) = (\D_{\f}  f, \wt g)$. Thus we conclude that $0\leq (\D_{\f}  f, \wt g)$ for all compactly supported functions $g$ with $g\geq 0$. \qed
\medskip

\medskip
For future reference we add a remark.

\Remark{\ZZ.12}  The operator $d^\f$ can be expressed in terms of the Hodge $d^*$-operator as
$$
d^\f f \ =\ d^*(f\f)
$$
and therefore
$$
dd^\f f\ =\ dd^*(f\f).
$$
To prove this, first note that if $v\in T_xX$ and $\a\in T^*_xX$ are metrically equivalent, then
$v\hk \f = *(\a\wedge *\f)$.  Hence, $d^\f f = \nabla f \hk \f = *(df \wedge *\f) = * (d(f*\f)-f(d*\f))
= *d*(f\f)-f*d*\f$, that is, 
$$
d^\f f\ =\ d^*(f\f)-f d^*\f
$$
so that the first equation holds if $\f$ is a harmonic form, and in particular if $\f$ is parallel.  Note also that for $\psi=*\f$
$$
d^\psi f\ =\ \pm *d(f\f) \and dd^\psi f \ =\ \pm *d^*d(f\f).
$$

\vskip .3in



\centerline{\bf  Appendix A.  Submanifolds which are $\f$-critical.}\medskip

\def\s{A}

This  appendix is not important for the remainder of the paper and can be skipped.
Here we establish a useful extension of Theorem \ZZ.3  to certain $\x$ which are not $\f$-planes.
Let $G \equiv  G(p,V)$ denote  the Grassmannian  of oriented p-planes in the inner product  space $V$, considered as the subset $G \subset \L_p V$ of unit simple vectors.

\Def{\ZZ.A.1}  Given $\f\in \L^pV^*$ an element $\x\in G$ is said to be a {\bf $\f$-critical
point} if $\x\in G$ is a critical point of the function $\f\bigl|_G$.
Equivalently, $\f$ must vanish on $T_\x  G \subset \L_pV$.
Let $G^{\rm cr}(\f)$ denote the set of $\f$ critical points.  \medskip

Note that if $\f$ is a calibration on $G$, i.e., $\sup \f\bigl|_G =1$, then
$$
G(\f)\ \subset \ G^{\rm cr}(\f)
$$
since for each $\x\in G(\f)$ the form $\f$ attains its maximum value 1 at $\x$.
Equation (\ZZ.9$)'$  extends from $G(\f)$ to $G^{\rm cr}(\f)$ as follows

\Prop{\ZZ.A.2}  {\sl   Suppose $\f\in \L^p V^*$  and   $\s \in\End(V)$.  Then for all $\x\in  G^{\rm cr}(\f)$}
$$
{\BM}_{\f} (\s)(\x)\ =\ (\tr_\x \s)\f(\x)
$$
\pf This is an immediate consequence of the more general Proposition \ZZ.A.4 below.\qed
\medskip

 Recall that at  a point $\x \in G$
there is  a canonical isomorphism:
$$
T_\x G\ \cong \ \Hom (\span \x, (\span \x)^\perp).
\eqno(\ZZ.A.1)
$$
On the other hand, $T_\x G$ is canonically a subspace of $\L_p V$.
It is exactly the subspace spanned by the first cousins of $\x$.
More specifically,  the isomorphism (\ZZ.A.1) associates to 
$L:\span \x \to (\span \x)^\perp$ the $p$-vector  $\Der_L \x$.

\Def{\ZZ.A.3} Let $\s\in\End (V)$ be a   linear map.  At each point $\x \in G$ we define
a tangent vector 
$$
\Der_{\wt A} \x\ \in\ T_\x G
$$
where  $\wt A =  P_{\x^{\perp}} \circ { \s }\circ P_\x$. This vector
field  $\x\to \Der_{\wt A} \x$  on $G$ is called the {\bf $A$-vector field}.
 
\Remark{} A straightforward calculation shows that this  $A$-vector field on $G$      
is the gradient vector field
$
 \n F_A
$
of the height  function $F_A : G \to \bbr$ given by $F_A(\x) =\langle A,
P_\x\rangle$.

\Prop{\ZZ.A.4}  {\sl Suppose $\f\in \L^p V^*$  and   $\s \in\End (V)$.
Then for all $p$-planes $\x\in G(p, V)$,}
$$
{\BM}_{\f} (\s)(\x)\ =\ (\tr_\x \s)\f(\x) + \f (\Der_{\wt A} \x)
\eqno(\ZZ.A.2)
$$
\pf
Pick an orthonormal basis for $\x$ and extend to an orthornormal basis of $V$.  
It suffices to prove (\ZZ.A.2) when   $A=b\otimes a$
with $a$ and $b$ elements of this basis. Using  formula (\ZZ.11) we see the following.

(1)\  If $a\in\x^\perp$, then all terms in (\ZZ.A.2) are zero.

(2)\  If $a\in\x$ and $b\in\x^\perp$, then $\wt A = A=b\otimes a$,   $\tr_\x A=0$,
and ${\BM}_{\f} (b\otimes a)(\x) = (a\wedge(b\hk\f)) = \f(b\wedge(a\hk\x)) = \f(\Der_A\x)$

(3)\  If  $a=b\in\x$, then ${\BM}_{\f} (\s)= \f(a\wedge(a\hk\x)=\f(\x)$ and $\tr_\x(A)=1$.
Since $\wt A=0$, equation (\ZZ.A.2) holds in this case.

(4)\  If $a,b \in\x$ and $a\perp b$, then $b\wedge(a\hk \x)=0$, and one sees easily that
all three terms in (\ZZ.A.2) are zero.
\qed

 \Remark{\ZZ.A.5} Proposition \ZZ.A.2  can be restated as
 $$
 \lambda_{\f}^*(\x)\ =\ \f(\x) P_\x \qquad {\rm for\ all\ }\ \ \x\in G^{\rm cr}(\f).
 \eqno{(\ZZ.A.3)}
 $$
 Conversely, if $ \lambda_{\f}^*(\x)\ =\ c P_\x$\  for some  $\x\in G^{\rm cr}(\f)$, then
 $c=\f(\x)$ and $\x$ is $\f$-critical.
 
 \pf  For all  $\x\in G(p,V)$ we have 
 $\langle  P_\x, \lambda_{\f}^*(\x)  \rangle =  ( \lambda_\f  P_\x)(\x)
  = (D_{P_{\x}^t}\f ) (\x)   =
  \f(D_{P_{\x}} \x)  = p\f(\x)$ 
since $D_{P_{\x}}\x=p\x$.  Therefore, $\lambda_{\f}^*(\x)=cP_\x$
 implies that $pc=p\f(\x)$ and equation (\ZZ.A.3) holds. 
 Equation (\ZZ.A.2) now implies that $\f(D_{\wt A}\x)=0$ for all $A\in \End(V)$ and, in particular,
  $\f(D_{L}\x)=0$ for all $L:\x\to \x^{\perp}$.  That is, $\f$ vanishes on $T_{\x}G\subset \L_p V$,
  i.e., $\x\in G^{\rm cr}(\f)$. \qed

\medskip

We now define an oriented submanifold $M$ of $X$ to be $\f$-{\bf critical} if $\oa T_x M \in G^{\rm cr}(\f)$ for all $x\in M$.  We leave it to the reader to use Proposition \ZZ.A.2  to establish the following extension of the previous results.

\Theorem{\ZZ.A.6} {\sl  Suppose $\f$ is a  $p$-form on a riemannian manifold 
$X$  and  $M\subset X$ is a $\f$-critical submanifold with mean curvature vector field $H$.     
Then for all $f\in C^\infty(X)$, 
$$
{\BM}_{\f} (\Hess f)\ =\ [\D_M (f) +H(f)]\f
$$
 when restricted to $M$. In particular, if $M$ is minimal, then on $M$}
 $$
 {\BM}_{\f} (\Hess f) \ =\ (\D_M f)\f
 $$
  
%
%
 %

\Ex{}  Let $\f={1\over 6}\{\o_I^2+\o_J^2+\o_K^2\}$ be the quaternion calibration on 
 $\bbh^n$. Then $\pm {1\over 3}$ are critical values and the  $\f$-critical submanifolds
 with critical value $\pm {1\over 3}$ include all complex Lagrangian submanifolds for any  complex
 structure defined by right multiplication by a unit imaginary quaternion (cf.  [U]).

   \vskip .3in

%
%

\centerline{\bf Appendix B.   Constructing $\phi$-plurisubharmonic  functions. }

\medskip
  
 Straightforward calculation shows that if $F(x)=g(u_1(x),...,u_m(x))$,  then
 $$
 d d^\phi F \ =\ \sum_{j=1}^m {\partial g \over \partial t_j} dd^\phi u_j
 + \sum_{i,j=1}^m {\partial^2 g \over \partial t_i\partial t_j} \BM_\phi(\nabla u_i \circ \nabla u_j)
 \eqno {(\ZZ.B.1)}
 $$

\Prop{\ZZ.B.1}  {\sl  If  $u_1,...,u_m$ are $\phi$-pluriharmonic  and $g(t_1,...,t_m)$ is convex,
then $F=g(u_1,...,u_m)$ is \fp. More generally, if ${\partial g \over \partial t_j} \geq 0$ for 
$j=1,...,m$ and $g$ is convex, then $F=g(u_1,...,u_m)$ is \fp whenever each $u_j$ is 
$\phi$-plurisubharmonic.
}
\pf Under our assumptions the  first term in equation (\ZZ.B.1) is $\geq 0$ on any $\xi\in G(\phi)$.
To show that the second term is $\geq 0$ is suffices to consider the case where the matrix
$(({\partial g \over \partial t_j}))$ is rank one, i.e., equal to $((x_ix_j))$ for some vector $x\in\bbr^n$.
Then the second term equals $\BM_\phi\{(\sum_i x_i\nabla u_i)\circ (\sum_j x_j\nabla u_j)\}$
which is $\geq0$  on $\xi\in G(\phi)$ by (\ZZ.13) and Corollary \ZZ.7. \qed
\medskip

We now analyze the case where $m=2$ and determine necessary and sufficient 
conditions for $F=g(u_1,u_2)$ to be $\phi$-plurisubharmonic.

\Lemma{\ZZ.B.2} {\sl  Fix $v,w\in \bbr^n$ and $\xi\in G(\phi)$.  Let $v_0$ and $ w_0$ denote the orthogonal
projections of $v$ and $w$ respectively onto   $\xi$ (considered as a $p$-plane in $\bbr^n$). Then}
$$
\BM_\phi(v\circ w)(\xi)\ =\   \langle v_0, w_0\rangle  .
$$
\pf
Write $v=v_0+v_1$ and $w=w_0+w_1$ with respect to the decomposition $\bbr^n = \span \xi
\oplus (\span \xi)^\perp$.  Then for $\xi \in G(\phi)$ we have
$$
\BM_\phi(v\circ w)(\x) \ =\  \phi\{(v_0+v_1)\wedge ((w_0+w_1) \hk \xi))\}
 \ =\   \phi\{(v_0+v_1)\wedge (w_0 \hk \xi)\}
\qquad \qquad $$ $$\qquad 
  \ =\   \phi(v_0\wedge (w_0 \hk \xi))  \ =\    \langle v_0, w_0\rangle  
\phi(\xi)  \ =\  
 \langle v_0, w_0\rangle  .
$$
where the third equality follows from the First Cousin Priinciple. \qed \medskip

By Lemma \ZZ.B.2 we have that for  $\xi\in G(\phi)$, 
$$\eqalign{
\BM_\phi\{a\,v\circ v +&2b \, v\circ w +c \,w\circ w\} (\xi)
\cr
&=\ a\|v_0\|^2 +2b \langle v_0, w_0\rangle   + c\|w_0\|^2
\cr
&=  \  \left\langle \left( \matrix{a &b \cr  b & c}\right),   \left( \matrix{\|v_0\|^2 & \langle v_0, w_0\rangle   
 \cr    \langle v_0, w_0\rangle   & \|w_0\|^2 }\right) \right\rangle
}
\eqno{(\ZZ.B.2)}$$

\Remark{\ZZ.B.3} A symmetric $n\times n$-matrix $A$ is $\geq 0$ iff $ \langle A, P\rangle  \geq 0$ for all rank-one symmetric
$n\times n$-matrices $P$.

\Remark{\ZZ.B.4}  The matrix  $ \left( \matrix{\|v_0\|^2 & \langle v_0, w_0\rangle     \cr    \langle v_0, w_0\rangle  & \|w_0\|^2 }\right)$
is rank-one iff  $v_0$ and $w_0$ are linearly dependent.

\Lemma{\ZZ.B.5}  {\sl  Let $v,w\in \bbr^n$  be linearly independent. Suppose that for every line
$$
\ell\ \subset \ \span\{v,w\}
$$
there exists a $(p-1)$-plane $\xi_0\subset \span\{v,w\}^\perp$ such that
$\ell \oplus\xi_0$ (when properly oriented) is a $\phi$-plane.  Then $\BM_\phi\{a v\circ v+2b v\circ w +cw\circ w \}$ is $\phi$-positive if and only if  
$\left(\matrix{a &b\cr b&c}\right) \geq 0$.
}

\pf
Necessity is already done.  For sufficiency fix $a,b,c$.  For each $\ell \subset  \span\{v,w\}$
let $\xi\in G(\phi)$ be the oriented $p$-plane $\ell\oplus \xi_0$ given in the hypothesis, and
 note that by equation (\ZZ.B.2)
 $$
\BM_\phi\{a v\circ v+2b v\circ w +cw\circ w \} (\xi)\ =\  \left\langle \left( \matrix{a &b \cr  b & c}\right),                \left( \matrix{v_{\ell}^2 &v_{\ell}w_{\ell}  \cr   v_{\ell}w_{\ell} & w_{\ell}^2 }\right) \right\rangle
\ \geq\ 0
 $$
where $v_{\ell}= \langle v,e\rangle$,  $w_{\ell}=  \langle w,e\rangle  e$, and $\ell = \span\{e\}$.  Now the matrix 
$\left( \matrix{v_{\ell}^2 &v_{\ell}w_{\ell}  \cr   v_{\ell}w_{\ell} & w_{\ell}^2 }\right) $
is rank-one, and every rank-one $2\times 2$ matrix, up to positive scalars, occurs in this family.
The result follows from Remark \ZZ.B.3. \qed

\Def{\ZZ.B.6}  A calibration $\phi$ on a manifold $X$ is called {\bf rich} (or {\bf 2-rich}) if for any 2-plane $P\subset  T_xX$ at any point x, and for any line $\ell\subset P$, there exists a $(p-1)$-plane
$\xi_0\subset P^\perp$ so that $\pm \ell\oplus \xi_0$ is a $\phi$-plane.

\Prop{\ZZ.B.7} {\sl   Let $(X,\phi)$ be a rich calibrated manifold.  Suppose   $u_1, u_2$ are 
$\phi$-pluriharmonic functions on $X$ 
with $\nabla u_1 \wedge \nabla u_2 \neq 0$ on a dense set. Then
for any $C^2$-function $g(t_1,t_2)$}
$$
F\ =\ g(u_1,u_2)  \in \fpsh \ \ \Leftrightarrow\ \ \ g {\rm \ \ is\ convex}
$$
\pf
Apply Proposition \ZZ.B.1, equation (\ZZ.B.2) and Lemma \ZZ.B.5.\qed

\Prop{\ZZ.B.8}  {\sl The Special Lagrangian calibration on a Calabi-Yau n-fold, $n\geq 3$, and the associative and coassociative calibrations on a $G_2$-manifold are rich calibrations.}

\pf
For the Special Lagrangian case it suffices to consider  $\phi ={\rm Re}(dz)$ on $\bbc^n$, $n\geq 3$.
Let $e_1, Je_1,... e_n, Je_n$ be the standard hermitian basis of $\bbc^n$.  By unitary invariance we may assume that $\ell = \span\{e_1\}$ and $P=\span\{e_1, \a Je_1+ \b e_2\}$. Then the $(p-1)$-plane 
$\xi_0 = -Je_2\wedge Je_3\wedge e_4\wedge\cdots \wedge e_n$ does the job.

Consider now the associative calibration $\phi(x,y,z) = \langle x\cdot y, z\rangle  $ on the imaginary octonians
${\rm Im}({\bbo}) = {\rm Im}({\bbh}) \oplus \bbh\cdot \epsilon$ where $\bbh$ denotes the quaternions.
By the transitivity of the group $G_2$ on $S^6=G_2/{\rm SU}(3)$ and the transitivity of SU(3)
on the tangent space, we may assume  $\ell = \span\{i\}$ and $P=\span\{i,j\}$ in Im$(\bbh)$.
We now choose $\xi_0 =  \epsilon\wedge ( i\cdot \epsilon)$.
For the coassociative calibration we choose $\xi_0 =k \wedge (i\epsilon)\wedge(k\epsilon)$
and note that $i\wedge \xi_0= i\wedge k \wedge (i\epsilon)\wedge(k\epsilon) $ is coassociative because its orthogonal complement   is $j\wedge\epsilon\wedge (j \epsilon) $ which is associative.
\qed

\medskip

We now give some examples and applications of the material above. We start with
Special Lagrangian geometry where the $\phi$-pluriharmonic functions are given by 
Proposition \AA.10.
Hence, we may apply Proposition \ZZ.B.7 to conclude the following. Let $u_1(z)$ and $u_2(z)$ be two
traceless hermitian quadratic forms on $\bbc^n$. (For example, take $u_1(z) =|z_1|^2-|z_2|^2$
and $u_2(z) =(n-2)|z_1|^2-|z_3|^2-\cdots-|z_n|^2$.)  Then $g(u_1(z),u_2(z))$ is \fp  if and only
if $g$ is convex.

\medskip
Formula (\ZZ.B.1) can be usefully appled to more general
functions $u_j$.   For example, in the Special Lagrangian case on
 $\bbc^n$ with $\phi={\rm Re}(dz)$, one has that $d d^\phi({1\over2}|z_k|^2) =\phi$, for any complex
 coordinate $z_k$  in any unitary coordinate system on $\bbc^n$. Hence  a linear combination of these functions have the property that $dd^\f u = c\f$ for some constant $c$.

\Prop{\ZZ.B.9} {\sl   Let $(X,\phi)$ be a rich calibrated manifold.  Suppose   $u_1, ..., u_n \in C^\infty(X)$  
 satisfy the equations $d d^\phi u_i= c_i \phi$ for constants $c_1,...,c_n$.  Then for any $C^2$-function $g(t_1,..., t_n)$}
$$
F\ =\ g(u_1,...,u_n)  \in \fpsh \ \ \Leftrightarrow\ \ \ \left\{\sum_{i=1}^nc_i{\partial g \over \partial t_i}\right\}
{\bf Id} \ +\   \left\langle \Hess_g ,  ((\langle (\nabla u_i)^\xi,  (\nabla u_j)^\xi  ))  \right\rangle             \ \geq \ 0
$$
for all $\phi$-planes $\xi$ at all points of $X$.

\def\s{\sigma}

\vfill\eject

\vskip .3in

\centerline{\bf  \BB. Convexity in Calibrated Geometries}\medskip

We suppose throughout this section that $(X,\f)$  is a non-compact, connected  calibrated manifold.
\medskip

\Def{\BB.1}  If $K$ is a compact subset of $X$, we define the {\bf $(X,\f)$-convex hull of } $K$ by
$$
\wh K\ \equiv\ \{x\in X : f(x) \leq \sup_K f \ \ {\rm for\ all }\ f \in  \cp\cs\ch(X,\f)\}
$$
If $\wh K = K$, then $K$ is called {\bf $(X,\f)$-convex.}

\Lemma{\BB.2}  {\sl Suppose $K$ is a compact subset of $X$.  Then $x\notin \wh K$ if and only if there exists a smooth non-negative  \fp function $f$ on $X$ which is identically zero on a neighborhood of $K$ and has $f(x)>0$.  Furthermore, if there exists a \fp function on $X$ which is strict at $x$, then  $f$ can be chosen to be strict at $x$.}

\pf  Suppose $x\notin \wh K$. Then there exists $g\in \PSH(X,\f)$ with $\sup_K g <0< g(x)$.
Pick $\varphi\in C^\infty(\bbr)$ with $\varphi\equiv 0$ on $(-\infty, 0]$ and with $\varphi>0$ and convex increasing on $(0,\infty)$.  Then $f=\varphi\circ g$ satisfies the required conditions.
Furthermore, assume $h\in \PSH(X,\f)$ is strict at $x$. Then take  $\ol g = g +\epsilon h$.
For small enough $\epsilon$,  $\sup_K \ol g <0<\ol g(x)$.  If $\varphi$ is also strictly increasing on $(0,\infty)$, then $f=\varphi \circ \ol g$ is strict at $x$.   \qed

\Theorem{\BB.3}  {\sl  The following two conditions are equivalent.
\smallskip
1)  \ \ If $K\subset\subset X$, then $\wh K \subset\subset X$.

\smallskip

2) \ \  There exists a \fp  proper exhaustion function $f$ on $X$.  }

\Def{\BB.4} If the equivalent conditions of Theorem \BB.3 are satisfied, then $(X,\f)$ is a 
{\bf convex calibrated manifold } and {\bf $X$ is \fc.}

\smallskip
\noindent
{\bf Proof that  $2) \Rightarrow 1)$:}
If $K$ is compact, then $c\equiv \sup_K f$ is finite and $\wh K$ is contained in the compact prelevel set $\{x\in X : f(x) \leq c\}$.

\medskip
\noindent
{\bf Proof that  $1) \Rightarrow 2)$:}  A \fp proper exhaustion function on $X$ is constructed as follows.
Choose an exhaustion of $X$ by compact  $(X,\f)$-convex subsets 
$K_1\subset K_2\subset K_3\subset \cdots $ with $K_m\subset K^0_{m+1}$ for all $m$.
By Lemma \BB.2 and the compactness of $K_{m+2}-K_{m+1}^0$, there exists a \fp function $f_m\geq0$ on $X$ with  $f_m$ identically zero on a neighborhood of $K_m$   and $f_m>0$ on $K_{m+2}-K_{m+1}^0$.
By rescaling we may assume $f_m>m$ on $K_{m+2}-K_{m+1}^0$. The locally finite sum 
$f=\sum_{m=1}^\infty f_m$ satisfies 2).\qed

\Lemma{\BB.5} {\sl Condition 2) in Theorem \BB.3 is equivalent to the {\rm a priori} weaker condition:
\smallskip

$2)'$  There exists a continuous proper exhaustion function $f$ on $X$ which is smooth and 

\quad \   \fp outside a compact subset of $X$.

\smallskip
\noindent
In fact if  $f$ satisfies $2)'$, then $f$ can be modified on a compact subset to be \fp on all of $X$.
Consequently, if $f$ satisfies $2)'$  and is strict outside a compact set, then its modification $\varphi\circ f$ is also strict outside a compact set.}

\pf    For large enough $c$, $f$ is smooth and  \fp outside the compact set    $\{x\in X : f(x)\leq c-1\}$.  Pick a convex increasing function $\varphi\in C^\infty(\bbr)$   with $\varphi \equiv c$ on a neighborhood of $(-\infty, c-1]$ and $\varphi(t)=t$ on $(c+1, \infty)$.  Then $\varphi\circ f$ is \fp on all of $X$ (in particular smooth) and equal to $f$ outside of the compact set  $\{x\in X : f(x)\leq c+1\}$.\qed

\Theorem{\BB.6} {\sl The following two conditions are equivalent:\smallskip

1)\ \ $K\subset\subset X  \ \Rightarrow \ 
\wh K\subset\subset X $, and $X$ carries a strictly \fp function.
\smallskip

2)\ \ There exists a strictly \fp proper exhaustion function for $X$.}

\Def{\BB.7} If the equivalent conditions of Theorem \BB.6 are satisfied, then   $(X,\f)$ is a
{\bf strictly convex calibrated manifold} or $X$ is {\bf strictly $\f$-convex}.
\medskip

\noindent
{\bf Proof of Theorem \BB.6.}  Suppose that $X$ is equipped with both a \fp proper exhaustion function
 $f$ and a strictly \fp function $g$. Then the sum $f+e^g$ is a strictly \fp exhaustion function.
 Now Theorem \BB.6 follows  immediately from Theorem \BB.3.\qed

\medskip

We shall construct many $\f$-convex manifolds in the course of our discussion
(See, in particular, \S \FF). However, we present some elementary examples here.

\Ex{1} Suppose $\f\in\L^p\bbr^n$ is a parallel calibration on $\bbr^n$.  Let
$f(x) = {1\over2}\|x\|^2$.  Then $dd^\f f=p\f$ and hence $f$ is a strictly \fp exhaustion.
That is,  $(\bbr^n, \f)$ is a strictly convex calibrated manifold.

\Ex{2}    Suppose $\f= dx_1\wedge\cdots\wedge dx_n$ on a domain $X$ in $\bbr^n$. 
Then $dd^\f f=(\D f)\f$ and $f$ is \fp if and only if $f$ is subharmonic.  Recall that if 
$K\subset\subset X$, then $\wh K = K \cup \{{\rm all\ the\  ``holes"\ in\ } K $ relative to $X\}$,
(connected components of $X-K$ which are relatively compact in $X$.
Thus $(X,\f)$  is strictly convex for any open set $X\subset \bbr^n$.

\medskip

It is instructive to extend this elementary example.

\Ex{3}      Suppose $\f= dx_1\wedge\cdots\wedge dx_p$ on a domain $X$ in $\bbr^n$. 
A function $f\in C^\infty(X)$ is  \fp if and only if $\D_x f\geq 0$ on $X$.
For a set $K\subset \bbr^n$, let $K_y$ denote the horizontal slice $\{x\in \bbr^p : (x,y)\in K\}$
of $K$.  Suppose that for each $y\in \bbr^{n-p}$, the horizontal slice $X_y$ has no
holes in $\bbr^p$.  Then $(X,\f)$ is strictly convex.  To prove this fact, it suffices to exhaust
$X$ by compact sets $K$ with the same property and show that each such $K$ 
is equal to its $(X,\f)$-hull.  Suppose $z_0=(x_0,y_0) \in X-K$.  Since $x_0$ is not in a hole
of $K_{y_0}$ in $\bbr^p$, we may choose (by Example 2) an entire subharmonic function
$g(x)$ with $g(x_0) >>0$ and $\sup_{K_{y_0}} g <<0$.  Now pick $\psi\in C^\infty_{\rm cpt}(\{y: |y-y_0|<\epsilon\})$ with $0\leq \psi\leq 1$ and $\psi(y_0)=1$.  Then $f(x,y) = g(x)\psi(y)$ is \fp 
and $f(z_0)=g(x_0)>>0$.  For $\epsilon$ sufficiently small, $\sup_K f\leq0$.  This proves  
$z_0$ does not belong to the $(X,\f)$-hull of $K$.

\Ex{4}  Let $\f=dx$ in $\bbr^2$ and set $X=\{(x,y) : x^2-c<y< x^2, \ |x|<1\}$.  Then $X$ is not $\f$-convex.  The closure of the hull of the compact subset
$K = ([-\e,\e]\times \{-\e\})\cup(\{\pm\e\}\times [-\e,0])$ of $X$ is easily seen to contain the origin.
Similarly, a domain of ``U''-shape, whose upper boundary along the bottom has a flat segment, is not $\f$-convex even though it is locally $\f$-convex
 (by Example 3).
\medskip

It is important to "weaken" this notion of strict convexity. 

\vfill\eject

\Theorem{\BB.8}  {\sl  The following two conditions are equivalent:\smallskip

 1)\ \ $K \subset \subset X  \ \Rightarrow \  \wh K\subset\subset X $, and there exists a strictly \fp function defined outside a compact subset of $X$  \smallskip

2)\ \ There exists a   \fp proper exhaustion function  on $X$ which is strict outside a compact subset
of $X$.}

\Def{\BB.9}  If the equivalent conditions of Theorem \BB.8 are satisfied,  then the calibrated manifold  $(X,\f)$ is  
{\bf strictly convex at  $\infty$} or $X$ is {\bf strictly $\f$-convex at $\infty$.}
 
 \Remark{}  This is not  the standard terminology used in complex geometry where such spaces are called ``strongly (pseudo) convex".

\medskip
\noindent
{\bf Proof of Theorem \BB.8.}  Obviously 2) implies 1).  We will prove that 1) implies the following weakening of 2).
\smallskip

$2)'$ There exists a continuous proper exhaustion function $f$ on $X$ which is smooth and strictly \fp outside a compact subset of $X$.

\smallskip
By  Lemma \BB.5, Condition $2)'$ implies Condition 2).

Now assume 1).  Since $K\subset \subset X $ implies $\wh K\subset\subset X $, we know from Theorem \BB.3 that there exists a \fp exhaustion function $f$ for $X$.  Let $g$ denote the strictly \fp function which is only defined outside of a compact set.  We can assume this compact set is $\{x\in X : f(x)\leq c\}$ for some large $c$.  Then $h\equiv \max\{f+e^g, c\}$ is a continuous proper exhaustion function which, outside the compact set $\{x\in X : f(x)\leq c\}$, is strictly \fp (in fact, equal to$f+e^g$).
This proves $2)'$ and completes the proof of the theorem.\qed

\Cor{\BB.10} {\sl $(X,\f)$ is strictly convex at $\infty$ if and only if Condition $2)'$ holds.}

\vskip .3in
\centerline{\bf Cores.}\medskip

In each non-compact calibrated manifold $(X,\f)$ there are certain distinguished subsets
which play an important role in the $\f$-geometry of the space.  (In complex manifolds
which are strongly pseudoconvex, these sets correspond to the compact exceptional
subvarieties.)  The remainder of this section is devoted to a discussion of these subsets.

Given a function $f\in\fpsh$, consider the open set
$$
S(f)\ \equiv\ \{x\in X : f \ \ {\rm is \  strictly \ }
\phi-{\rm plurisubharmonic  \ at\ \  } x \}
$$
 and the closed set
$$
W(f)\ \equiv \ X-S(f).
$$
Note that 
$$
W(\lambda f + \mu g) \ =\ W(f) \cap W(g) 
$$
for $f,g \in \fpsh$ and $\lambda, \mu >0$.

\Def{\BB.11}  The {\bf core} of $X$ is defined to be the intersection
$$
\Core(X) \ \equiv\ \bigcap W(f)
$$
over all $f\in\fpsh$. The {\bf inner core} of $X$ is defined to be the set ${\rm InnerCore}(X)$
of points $x$ for which there exists $y\neq x$ with the property that $f(x)=f(y)$
for all $f\in\fpsh$.

\Prop{\BB.12}  \qquad InnerCore$(C) \ \subset \ \Core(X)$.

\pf If $x\notin \Core(X)$, then there exists $g\in \PSH(X,\f)$ with $g$ strict at $x$.
Suppose $y\neq x$. Then if $\psi$ is compactly supported in a small neighborhood
of $x$ missing $y$, and $\psi$ has sufficiently small second derivatives, one has 
$f=g+\psi\in \PSH(X,\f)$. Obviously for such $f$, the values $f(x)$ and $f(y)$ can be
made to differ, so therefore $x\notin $ InnerCore$(X)$. \qed

\Prop{\BB.13}  {\sl  Every compact $\f$-submanifold $M\subset X$ is contained in the inner core.}

\pf  Each $f\in\PSH(X,\f)$ is subharmonic on $M$ by Theorem \AA.4.  
Hence,   $f$ is constant on $M$. \qed

\Prop{\BB.14}  {\sl  Suppose $X$ is $\f$-convex.  Then $\Core(X)$ is compact
if and only if $X$ is strictly $\f$-convex at $\infty$,
and $\Core(X)=\emptyset$ if and only if $X$ is strictly \fc.}

\pf
If  $X$ is strictly $\f$-convex at $\infty$, then choosing $f$ to satisfy 2) in Theorem \BB.8, we see that the Core$(X) \subset W(f)$ is compact.  Obviously, strict $\f$-convexity implies that $\Core(X)=\emptyset$.

Conversely, if $\Core(X)$ is compact, then in the construction of the \fp exhaustion function in the proof of Theorem \BB.3 we may choose
$$
 K_1\ =\ \wh{\Core(X)}
 $$
Then by the definition of $\Core(X)$ and Lemma \BB.2, each of the functions $f_m$ in that proof can be chosen to be strictly \fp 
on $K_{m+2}-K^0_{m+1}$. Hence the exhaustion $f=\sum_m f_m$ is strictly \fp outside a compact 
set containing the core.\qed
\medskip

A slight modification of this construction gives the following general result.

\Prop{\BB.15} {\sl  Suppose  $X$ is strictly $\f$-convex at $\infty$, and $K\subset X$ is
a compact, $\f$-convex subset containing the core of $X$.  Let $U$ be any
neighborhood of $K$. Then  there exists a proper \fp exhaustion function $f:X\to \bbr^+$
 which is strictly \fp  on $X-U$, and  identically zero on a neighborhood of $K$.
}
\pf    
Choose $K_1=K$ in the construction of the \fp exhaustion function given in the proof of Theorem \BB.3.
Let $K_\epsilon$ denote the compact $\e$-neighborhood of $K$. Then
$$
K\ =\ \bigcap_{\e>0} \wh{K}_\e.
\eqno{(\BB.1)}
$$
If $x\in \bigcap_{\e>0} \wh{K}_\e$, then for each $f\in \PSH(X,\f)$, we have $f(x) \leq \sup_{K_\e} f$.
However,   $\inf_{\e} \sup_{K_\e} f = \sup_K f$, and we conclude that $x\in \wh K$.  Thus we can choose $K_2 \equiv  \wh K_\e$ in our construction of $f$, and for  small enough $\e$ we have $K_2\subset U$ as well as $K_1\subset K_2^0$. The proof is now completed as in the proof of Proposition \BB.14.  \qed
\medskip

Obviously, many question concerning
$$
{\rm InnerCore}(X)\ \subset\ \Core(C)\ \subset\ \wh{\Core(X)}
$$
remain to be answered.

\vfill\eject


\centerline{\bf Appendix \B. Structure of the Core.}\medskip

Let $(X,\f)$ be a calibrated manifold and consider the set
$$
\cn \ \equiv\ \{ \x \in G(\f) : (\ch^\f f)(\x)\ =\ 0 \ \ {\rm for\ all\ }
f\in\fpsh\}.
$$
\Prop{\BB.\B.1} {\sl  Let $\pi:G(\f)\to X$ denote the projection. Then}
$$
\pi(\cn)\ =\ \Core(X).
$$
\pf
Suppose $x\notin\Core(X)$.  Then by definition there exists $f\in\fpsh$ 
with $(\ch^\f  f)(\x)>0$ for all $\x\in  \pi^{-1}(x)$.  Hence, $x\notin
\pi(\cn)$.

Conversely, suppose $x\notin \pi(\cn)$. Then for each $\x\in \pi^{-1}(x)$
there exists $f_\x\in\fpsh$ with $f(\x)>0$. Let $W_\x=\{\eta\in 
\pi^{-1}(x): f_\x(\eta)>0\}$ and choose a finite cover $W_{\x_1},..., 
W_{\x_\ell}$ of $\pi^{-1}(x)$. Then $f\equiv f_{\x_1}+\cdots+f_{\x_\ell}$
is strictly \fp at x, and so $x\notin \Core(X)$.  \qed

\Prop{\BB.\B.2} {\sl If $\x\in \cn$, then for each vector $v\in\span \x$,}
$$
df(v) \ =\ 0 \qquad{\sl  for\ all }\ \  f\in\fpsh
\eqno(\BB.\B.1)$$
\pf
Suppose $f\in\fpsh$ and set $F=e^f$.  Then $F\in \fpsh$, and   by  
 equation (\ZZ.14) and Corollary \ZZ.7 we see that
$
0\ =\ (\ch^\f F)(\x)\ =\ e^f\{df\w d^\f f + \ch^\f f\}(\x)\ 
=\ e^f\{df\w d^\f f \} (\x) \ = \   
e^f|\n f\hk\x|^2
$.  \qed
\medskip

\Def{\BB.\B.3}  The {\bf tangential core} of $X$ is the set
$$
T\Core(X)\ \equiv\ \{ v\in TX: v \neq0\  {\rm and\ satisfies\ condition \
(\BB.\B.1)} \} $$
Thus $T\Core(X) \subset TX$ is a subset defined by the vanishing of 
the family of smooth functions $df:TX\to \bbr$ for $f\in\fpsh$.
Propositions \BB.\B.1 and \BB.\B.2 show that the restriction of the bundle map
$p:TX\to X$ gives a surjective mapping
$$
p:T\Core(X) \to \Core(X)
$$
and for each $x\in X$, The vector space $T_x\Core(C)\equiv p^{-1}(x)$
contains the non-empty space generated by all $v\in \span \x$ for $\x\in
\cn_x$.

Consider a point $v\in T\Core(X)$ and suppose we have functions
$f_1,...,f_\ell\in\fpsh$ such that $\n df_1,..., \n df_\ell$ are linearly
independent at $v$. Then $T\Core(C)$ is locally contained in the
codimension-$\ell$ submanifold $\{df_1=\cdots=df_\ell=0\}$. If additionally
we assume that $df_1,...,df_\ell$ are linearly independent vectors at 
$x=p(v)$, then $\Core(X)$ is locally contained in the subset
$\{f_1={\rm constant}\}\cap\cdots\cap\{f_\ell={\rm constant}\}$.
NO

\vfill\eject


\centerline{\bf  Appendix B.  Examples of Complete  Convex Manifolds and Cores}

\medskip

In \S \FF \ (Theorem \FF.4) we shall show that there are many strictly $\f$-convex domains in any calibrated manifold
$(X,\f)$.  They can have quite arbitrary topological type within the strictures imposed by Morse Theory
and $\f$-positivity of the Hessian.  However, it is  also interesting geometrically to ask for convex manifolds which are complete.
 
 In fact, there exist enormous families of  complete calibrated manifolds $(X,\phi)$ with $\nabla\phi=0$ 
 which are  strictly $\phi$-convex at infinity.  For example any asymptotically locally euclidean  (ALE)
 manifold with SU(n), Sp(n), G$_2$, or Spin$_7$ holonomy is such a creature, since the 
 radial function on the asymptotic chart at infinity is strictly  convex.  For the general
 construction of such spaces the reader is referred to the book of Joyce [J].
 
 However, some manifolds of this type have been explicitly constructed, and in these cases
 one can explicitly construct \fp exhaustion functions and identify the cores. We indicate how
 to do this below.
 
 We begin however with an observation in dimension 4. Every crepant resolution of singularities
 of $\bbc^2/\G$ admits Ricci-flat ALE K\"ahler metric.  On each such manifold there exists an $S^2$-family of parallel calibrations
$$
 \cc\ =\  \{u\omega + v\vf+\ w \psi : u^2+v^2+w^2=1\}
 $$
where $\o$ is the given K\"ahler form, 
 $\varphi= {\rm Re}\{\Phi\}$  and $\psi  ={\rm Im}\{\Phi\}$ and $\Phi$ is a parallel section of the canonical bundle $\kappa_X$.   Let $E=\pi^{-1}(0)$ be the exceptional locus of the resolution. Then for any $\phi\in \cc$  we have
 $$
 \Core(X,\phi) \ =\ \cases {E & if $\phi=\o$ \cr \emptyset & otherwise}
 $$
 This follows from the fact that each $\phi\in \cc$ is in fact the K\"ahler form for a complex structure
 on $X$ compatible with the given metric.  With this complex structure $X$ is pseudo-convex, and by the Stein Reduction Theorem  (cf.   [GR, p. 221])  we know its core is the union of its compact complex subvarieties.  For $\phi\neq \o$ there are no such subvarieties since by the Wirtinger inequality  (cf.  [L$_{1,2}$]),  applied to $\phi$, they would necessarily be homologically mass-minimizing, and by the same result applied to $\o$  any such subvariety is $\o$-complex (and therefore a  component of $E$).
 
 \Ex {1. (Calabi Spaces)}  Let $X\to \bbc^n/\bbz_n$ be a crepant resolution of $\bbc^n/\bbz_n$
 where the action on $\bbc^n$ is generated by scalar multiplication by $\tau = e^{2\pi i/n}$.
Following Calabi [C] we define the function $F:\bbc^n/\bbz_n\to \bbr$ by
$$
F(\rho) \ =\ \root n \of {\rho^n+1} + {1\over n}\sum_{k=0}^{n-1} \tau^k \log\left(\root n \of {\rho^n+1} -\tau^k\right)
$$
where  $\rho \equiv \|z\|^2$ (pushed down to $\bbc^n/\bbz_n$), and the log is defined by 
choosing $\arg \zeta \in (-\pi,\pi)$.
 We then define a Kahler metric on $\bbc^n/\bbz_n-\{0\}$ by setting
 $$
 \o \ =\ {1\over 4} d d^c F
 $$
 Calabi shows that {\sl this metric is Ricci-flat  and (when pulled back)  extends to a Ricci flat metric on $X$.  The parallel form $\Phi = dz_1\wedge\cdots\wedge dz_n$ extends to a parallel section of } $\kappa_X$.   This metric is given explicitly on $\bbr^{2n}/\bbz_n$ by
 $$
ds^2\ =\  F'(\rho) |dx|^2 + \rho F''(\rho) dr \circ d^c  r
$$
 where $r=\|x\|$.  Define $G(\rho)$ by settiing $G'(\rho)=F'(\rho)+\rho F''(\rho)$ and $G(0)=0$.
 Then direct calculation shows that
 $$
 dd^\f G \ =\ 2n \phi
 $$
 where $\phi = {\rm Re} \{\Phi\}$.  Hence, {\sl  $X$ is a complete,  strictly $\phi$-convex manifold.}
 
 \Ex{2. (Bryant-Salamon Spaces)}  Let $P$ denote the principal Spin$_3$-bundle of $S^3$ and 
 $$
 S\ \equiv\ P \times_{{\rm Sp}_1} \bbh
 $$
 the associated spinor bundle, where $\bbh $ denotes the quaternions.
Bryant and Salamon have explicitly constructed  a complete riemannian metric with
G$_2$-holonomy on the total space of $S$. (See [BS, page 838, Case ii].)  Let $\rho=|a|$
for $a\in \bbh$ (pushed-down to $S$) and let $Z\subset S$ denote the zero section.  Then a direct calculation shows that the function
$$
F(\rho) \ =\ (1+\rho)^{{5\over 6}} \ \ {\sl  is\ strictly\ } \varphi-{\sl plurisubharmonic\ on \ } S-Z
$$
 where $\varphi$ denotes the associative calibration on $S$.   Since $Z$ is an associative  submanifold
 we conclude that
 $$
 \Core(S)\ =\ Z
 $$
 In an analogous fashion the authors construct a complete riemannian metric with 
 Spin$_7$-holonomy on the total space $\wt S$ of a spinor bundle over $S^4$.
 (See  [BS, page 847, Case ii].)  A similar calculation shows that there exists an exhaustion
 function which is strictly $\Phi$-plurisubharmonic on  $\wt S - \wt Z$ where $\Phi$ denotes the Cayley
 calibration  $\wt Z$ the zero-section of $\wt S$. Since $\wt Z$ is a Cayley submanifold, we conclude that $$
 \Core (\wt S)\ =\ \wt Z
 $$

 \vfill\eject


\centerline{\bf \CC. Boundary Convexity.}
 
\medskip

Suppose $\O\subset\subset X$ is an open set with smooth boundary $\bo$, where $(X,\f)$ is a            non-compact calibrated manifold.  A $p$-plane $\x\in G(\f)$ at a point  $x\in \bo$ will be called  {\bf tangential} if $\span\x\subset T_x\bo$.

  \Def{\CC.1}  Suppose that $\rho$ is a {\bf defining function} for  $\bo$,
  that is, $\rho$ is a smooth function defined on a neighborhood 
  of $\overline{\Omega}$ with $\Omega =\{x: \rho(x)<0\}$ and  $\nabla \rho \neq 0$ on $\partial \Omega$.
  If
  $$
  dd^\phi\rho(\xi) \ \geq\ 0 \ \ \  {\rm for\ all\  tangential\ } \xi\in G_x(\phi), \ x\in \bo,
   \eqno{(\CC.1)} 
 $$
  then  $\partial \Omega$ is  called {\bf $\phi$-convex}.
  If the inequality in (\CC.1) is strict for all $\xi$, then 
  $\partial \Omega$ is  called {\bf strictly $\phi$-convex}.  If 
  $ dd^\phi(\xi) =0$ for all $\xi$ as in (\CC.1), then    $\partial \Omega$ is {\bf $\phi$-flat}.
  
  Each of these conditions is a local condition on $\bo$.  In fact:

  \Lemma {\CC.2} {\sl  Each of the three conditions in Definition \CC.1  is independent of the choice of defining function $\rho$.
  In fact, if $\overline{\rho}=u\rho$ is another choice with $u>0$ on $\partial \Omega$, then on $\bo$ }
  $$
  dd^\phi(\overline{\rho})(\xi)\ =\ udd^\phi(\rho)(\xi) \ \ \ {\rm for \ all \  tangential\ }\ \ \xi\in G(\phi)
    \eqno{(\CC.2)}
    $$
    
  \pf
  Note that $  dd^\phi(\overline{\rho})(\xi)=\ udd^\phi(\rho)(\xi) +  \rho dd^\phi(u)(\xi)
  + \BM_\phi(\nabla u\circ \nabla \rho)(\xi)$. Now the middle term  drops because $\rho=0$
  on $\bo$.  Furthermore, $\BM_\phi(\nabla u\circ \nabla \rho)(\xi) = 
  {1\over 2}\nabla \rho\wedge(\nabla u \hk \phi)(\xi) +
  {1\over 2}\nabla u\wedge(\nabla \rho \hk \phi)(\xi)
  = {1\over 2}\phi(\nabla u\wedge(\nabla \rho \hk \xi)) +
   {1\over 2}\phi(\nabla \rho\wedge(\nabla u \hk \xi))$. Now $\nabla \rho\perp \span (\xi)$ implies 
   $\nabla \rho \hk \xi =0$ and, by the First Cousin Principle, that $\phi$ vanishes on $\nabla \rho\wedge(\nabla u\hk \xi)$.\qed

  \Cor{\CC.3} {\sl  Assume $\f\in \L^p\bbr^n$ is a calibration.  Suppose $\bo$ is (strictly) $\phi$-convex in $\bbr^n$, and locally near a point $p\in \bo$, 
  let   $\bo$ be graphed over its tangent space by a function $x_n=u(x')$ for linear coordinates
  $(x',x_n)$ on $\bbr^n$.  Then each nearby hypersurface: $x_n = u(x')+c$ is also (strictly) $\phi$-convex.
  }\medskip
    The next lemma will be used to establish the main result of this section.
  
  \Lemma {\CC.4} {\sl Suppose $\rho$ is a smooth real-valued function defined near a point $x$
  in a calibrated manifold $(X,\f)$, and that $\psi:\bbr\to\bbr$ is smooth near $\rho(x)$.
  Then:
  $$
  \BM_{\f}(\Hess_x\psi(\rho))(\x)\ =\ \psi'(\rho)\BM_\f(\Hess_x\rho)(\x) + |\nabla\rho|^2\psi''(\rho) \cos^2\theta(\x)
  \eqno{(\CC.3)}
  $$
  where  $\cos^2\theta(\x) =\langle P_{\span \nabla \rho}, P_{\span \x} \rangle$.  }
  
  \pf
  First note that 
  $$
  \Hess \psi(\rho)\ =\ \psi'(\rho)\Hess \rho + \psi''(\rho)\nabla\rho\circ \nabla\rho, 
  \eqno{( \CC.4)}
  $$
  and then set $n=\nabla \rho/|\nabla \rho|$ so that 
  $
  \nabla\rho\circ \nabla\rho = |\nabla\rho|^2 n\circ n = |\nabla\rho|^2 P_{\span \nabla \rho}.
  $
  For each $\x\in G_x(\f)$, taking the inner product of (\CC.4) with $P_{\span\x}$ (orthogonal projection onto $\span \x$) yields (\CC.3) because of Lemma \AA.15b.\qed
  \medskip
  
  We now come to the main result of this section.

  \Theorem{\CC.5}   {\sl   Let $\Omega\subset\subset X$ be a compact domain with strictly \fc boundary.          Suppose $\delta = -\rho$ is an arbitrary ``distance function'' for $\bo$, i.e.,   $\rho$
  is an arbitrary  defining function for $\bo$.
  Then $-\log \,\delta$ is strictly \fp   outside a compact     subset of $\O$. Thus, in particular, the domain $\Omega$ is strictly $\phi$-convex at $\infty$.   }
  $ \over $
  
 \pf   
 Set $\psi(t)=-\log(-t)$ for $t<0$.  Note that $\psi'(t)=-1/t $ and  $\psi''(t)=-1/t ^2$, so that $\psi'(\rho)=
 1/\delta$  and $\psi''(\rho)= 1/\delta^2$.  Consequently, by Lemma 3.4, at each point $x\in\O$
 near $\bo$, we have
  $$
  dd^\phi (-\log\, \delta)(\x)\ =\ {1\over {\delta}} dd^\phi(\rho)(\x) + {|\nabla\rho|^2\over {\delta^2}} \cos^2 \theta  (\x)
  \eqno{(\CC.4)}
  $$ 
 for all $\x\in G(\f)$.  Note that at $x\in \bo$,  $\cos^2\theta(\x) = 
 \langle P_{\span \nabla \rho}, P_{\span \x} \rangle$ vanishes if and only if $\x$ is tangential to $\bo$. Consequently, 
the inequality $|\cos\theta|<\epsilon$ defines a fundamental neighborhood system for $G(p,T\bo)\subset G(p,TX)$.  By restriction     $|\cos\theta|<\epsilon$ defines a fundamental neighborhood system for $G(\phi)\cap G(p,T\bo)\subset G(\phi)$.  The hypothesis of strict  $\phi$-convexity for $\bo$ implies that there exists $\overline \epsilon >0$ so that
  $(dd^\phi \rho)(\x) \geq \overline\epsilon$ for all $\phi$-planes $\x$  at points of $\bo$ with $|\cos\theta|<\epsilon$
  for some $\epsilon>0$. (Note that if there are no $\phi$-planes tangent to $\bo$ at a point
  $x$, then there are no $\phi$-planes with $|\cos\theta|<\epsilon$ for sufficiently small
  $\epsilon$ in a neighborhood of $x$.)   Consequently, we have by equation (\CC.4)  that 
   $$
  dd^\phi (-\log\delta)(\x) \ \geq \ {\overline \epsilon \over 2\delta}
  $$
    near $\bo$ for all  $\phi$-planes $\x$  with $|\cos\theta|<\epsilon$, where $\theta$ is defined 
  as above with $\bo$ replaced by the nearby   level sets  of $\rho$.

  Now choose $M>>0$ so that $dd^\phi(\rho)(\x)\geq -M$ in a neighborhood of $\bo$ for all 
  $\x$.  Then, by (\CC.4) 
 $$
   dd^\phi (-\log\, \delta)(\x) \ \geq\   -{M\over \d}+  {1\over \d^2}  |\nabla \rho|^2 \cos^2\theta.
  $$
  If $|\cos\theta|\geq \epsilon$, this 
   is positive in a neighborhood of $\bo$ in $\Omega$.  This proves that $-\log \d$ is strictly \fp near $\bo$.    By Corollary \BB.10 the domain $\O$ is strongly $\f$-convex.
  \qed  
 
 \medskip
 
 Although a defining function for a strictly \fc boundary may not be $\f$-plurisub-harmonic, for some applications the following may prove useful.
 
 \Prop{\CC.6} {\sl  Suppose $\O\subset\subset X$ has strictly \fc boundary $\bo$ with defining function     $\rho$.  Then, for $A$ sufficiently large, the function $\overline \rho \equiv \rho + A \rho^2$ is strictly \fc in a neighborhood of $\bo$  and also a defining function for $\bo$.}

 \pf 
 By Lemma \CC.4
 $$
 \BM_\f(\Hess \overline\rho)(\x)\ =\ (1+2A\rho) \BM_\f(\Hess \rho)(\x) +2|\nabla \rho|^2A\cos^2\theta(\x)
 \qquad{\rm for\ all\ } \x\in G(\f)
 \eqno{(\CC.5)}
 $$
 where $\cos^2\theta(\x) = \langle P_{\span \nabla \rho}, P_{\span \x} \rangle$.  As noted in the proof 
 of Theorem \CC.5, there exist $\e, \overline\e>0$ so that $\BM_\f(\Hess \rho)(\x)\geq \overline \e$ if
 $|\cos\theta(\x)|<\e$, because of the strict boundary convexity.  Therefore $\BM_\f(\Hess \overline\rho)(\x) \geq (1+2A\rho)\overline \e$ if $\x\in G(\f)$ with $|\cos\theta(\x)|<\e$.  Choose a lower bound $-M$ for $\BM_\f(\Hess \rho)(\x)$ over all $\x\in G(\f)$ for a neighborhood of $\bo$.  Then by (\CC.5), 
 $\BM_\f(\Hess \rho)(\x)\geq -(1+2A\rho)M + 2|\nabla\rho|^2 A \e^2$ for $\x\in G(\f)$ with 
 $|\cos \theta(\x)|\geq \e$.  For $A$ sufficiently large, the right hand side is $>0$ in some neighborhood of $\bo$.\qed
 \medskip

One might hope for a converse to Theorem \CC.5, e.g.,   if the domain $\O$ is \fc then the  boundary
is \fc.  However, elementary examples show that this is false. 

\Ex{}  Let $\f\equiv dx\wedge dy$ in $\bbr^3$ as in Example 3 of section \BB.  Let $X$ denote the solid torus  obtained by rotating the disk $\{(y,z) : y^2+(z-R)^2 < r^2\}$ about the $y$-axis.  Since each slice $X_z$ has no holes in $\bbr^3$, the domain $X$ is $\f$-convex
(cf. Example 3 of \S \BB).  However, the boundary torus $\partial X$ is \fc  if and only if  $2r\leq R$.  This follows from  an elementary calculation which uses the obvious defining function and   Definition  \CC.1  (or by using Proposition \CC.12 below)

\Qu{\CC.7}  For which strictly convex calibrated manifolds is it true that \fc subdomains have \fc boundaries?  More generally, when is the $\f$-convexity of a domain a local condition at the 
boundary?
  \medskip
  
  A weak partial converse to Theorem \CC.4 is given by the following.
  
  \Prop{\CC.8}  {\sl  Suppose the calibration is parallel, and set $\d=\dist(\bullet, \bo)$.
  If $-\log \d$ is strictly \fp near $\bo$, then $\bo$ is \fc.
  }

\Note{\CC.9}  Examples show that the strict convexity of $-\log\d$ near $\bo$ is stronger than 
$\f$-convexity for $\bo$.

\pf   Set $\rho =-\delta$ on $\overline\O$   near  $\bo$.
  Suppose that $\bo$ is not $\phi$-convex. Then there exist $x\in \bo$ and 
  $\x_x\in G_x(\phi)$ with $\span(\x)\subset  T_x(\bo)$ and $(dd^\phi \rho)(\x_x)<0$.
  Let $\gamma$ denote the geodesic segment in $\Omega$ which emanates orthogonally from $\bo$ at $x$. Since $\delta$ is the distance function, $\gamma$ is an integral curve of $\nabla \delta$.  Let 
  $\x_y$, $y\in \gamma$ denote the parallel translation of $\x_x$ along $\gamma$.  Then $\x_y$ is a 
  $\phi$-plane with $\span(\x_y)\perp \nabla \rho$ for all $y$.  By formula (\CC.4), siince 
  $\cos\theta(\x_y)=0$, we have
  $$
  dd^\phi (-\log\, \delta)(\x_y) \ =\ {1\over \delta} dd^\phi(\rho)(\x_y) \ <\ 0
  $$
  for all $y$ sufficiently close to $x$.  Hence, $-\log \delta$ is not \fp near $\bo$.  \qed
  
\medskip

  The $\phi$-convexity of a boundary can be equivalently defined in terms of its second fundamental form.  Note that if $M\subset X$ is a smooth hypersurface with  a chosen unit normal field $n$ we have a quadratic form $II$ defined on $TM$ by 
  $$
  II(V,W) \ =\ \langle  B_{V,W}, n \rangle
  $$  
  where $B$ denotes  the second fundamental form of $M$ discussed in \S \AA.
  For example, when $H=S^{n-1}(r)\subset \bbr^n$ is the euclidean sphere of radius $r$, oriented by the outward-pointing unit normal,  we find that $II(V,W) = -{1\over r} \langle V,W\rangle$.

For the sake of completeness we include a proof of the following standard fact.

\Lemma{\CC.10} {\sl Suppose $\rho$ is a defining function for $\O$ and let II denote the second
fundamental form of the   hypersurface $\bo$ oriented by the outward-pointing normal.
Then}
$$
\Hess\, \rho\,\bigl|_{T\bo} \ =\ -|\nabla \rho| \, II
$$
\pf  Suppose $e$ is a tangent field on $\bo$. Extend $e$ to a vector field tangent to the level sets
of $\rho$. By definition $II(e,e)=\langle \nabla_e e, n\rangle$ where 
$n = \nabla \rho/|\nabla \rho|$ is the outward normal.
Then $(\Hess\, \rho)(e,e) = e(e \rho) -(\nabla_e e)\rho = -(\nabla_e e)\rho = -\langle \nabla_e e, \nabla \rho\rangle = - |\nabla \rho| \langle \nabla_ee, n\rangle$.
\qed

  \Remark{}
    Recall that a defining function $\rho$ for $\Omega$ satisfies  $\|\nabla \rho\| \equiv 1$ in a neighborhood of $\bo$ if and only if $\rho$ is the signed distance to $\bo$ ($<0$ in $ \Omega$ and $>0$ outside of $\Omega$).   In fact any function $\rho$ with $\|\nabla \rho\| \equiv 1$  in a riemannian manifold is, up to an additive constant,   the distance function to (any) one of its level sets.
  In this case it is easy to see that
      $$
    \Hess_\rho \ =\ \left( \matrix{0 & 0 \cr 0 & - II}       \right)
   \eqno{(\CC.6)} 
   $$
   where $II$ denotes the second fundamental form of the hypersurface $H=\{\rho =\rho(x)\}$
   with respect to the normal $n= \nabla \rho$ and the blocking in (3) is with respect to the 
   splitting $T_xX = \span (n_x) \oplus T_x H$.     For example let $\rho (x) =\|x\|\equiv r$ in $\bbr^n$. Then direct calculation shows that $\Hess_r 
 = {1\over r}(I-\hat x \circ \hat x)$ where $\hat x = x/r$.

  \Cor{\CC.11}  {\sl   For all tangential $\x\in G(\f)$}
  $$
  (dd^\f \rho)(\x) \ =\ -\|\nabla \rho\| \tr_\x II.
  $$
  \pf Apply Lemma \AA.10.\qed
  \medskip
  
  As an immediate consequence we have
  
  \Prop{\CC.12}  {\sl  Let $\Omega\subset X$ be a domain with smooth boundary $\bo$ oriented by
  the outward-pointing normal. Then $\bo$ is \fc if and only if its second fundamental form satisfies
  $$
  \tr_\x II \ \leq \ 0 
  $$
  for all $\phi$-planes $\xi$ which are tangent to $\bo$.  This can be expressed more geometrically by saying that
  $$
   \tr\left\{ B \bigr|_\xi\right\} \ \ \ {\rm must\ be\  inward-pointing\  }
  $$
  for all tangential $\phi$-planes $\xi$.
  }
    \medskip

\vfill\eject


\vskip .3in

\centerline{\bf \DD. Positive Currents in Calibrated Geometries.}

\medskip

The important classical notion of  a positive current on a complex manifold has an analogue
on any calibrated manifold.  This concept was introduced   in section II
of [HL$_3$].  We begin this Section by reviewing that  material with some of the terminology and notation updated.  

On a calibrated manifold $(X,\f)$ we have:
\smallskip

\qquad  a)\ \ $\f$-submanifolds,\smallskip

\qquad  b)\ \ rectifiable $\f$-currents, and\smallskip

\qquad   c)\ \ $\f$-positive (or $\lp$-positive)  currents.\smallskip

A $\f$-submanifold is, of course, a smooth oriented submanifold $M$ whose oriented
tangent space is a $\f$-plane at every point, i.e., $\oa M_x \equiv \oa T_x M\in G(\f)$ for all $x\in M$.

Suppose $T$ is a locally rectifiable $p$-dimensional current ([F]) on $X$. Then its generalized tangent space  is a unit simple vector $\oa T\in G(p, TX)$  
at $\|T\|$ almost every point, where $\|T\|$ denotes the generalized 
volume measure associated with $T$.

\Def{\DD.1}  A {\bf  rectifiable $\f$-current}  is a locally recitifiable current
$T$ with $\oa T \in G(\f)$ for $\|T\|$- a.a. points in $X$.  A {\bf  $\f$-cycle}
is a rectifiable $\f$-current which is $d$-closed.

\medskip

\Remark{\DD.2} We shall see below (Theorem \DD.9) that $\f$-cycles always have a particularly nice
local  structure.  The  strongest  result of this kind occurs in the K\"ahler case (where $\f=\o^p/p!$)
where  a theorem of  J. King [K] states that each 
$\f$-cycle is a positive holomorphic cycle, i.e., a locally finite sum of $p$-dimensional complex analytic subvarieties  with positive integer coefficients.   On a general calibrated
manifold $(X,\f)$ one can also consider $d$-closed recitifiable currents $T$ with $\pm \overrightarrow T\in G(\f)$ for $\|T\|$-a.a. points.  In the K\"ahler case $T$ must be a holomorphic chain by a theorem
of Harvey-Shiffman [HS], [S].  However, nothing is known about the structure of such currents for any of the other standard calibrations.

\medskip

An understanding of the definition of a $\f$-positive current is a little more complicated.

Recall(Federer [F]) that a current $T$ is {\bf representable by integration} if $T$ has measure coefficients when expressed as a generalized differential form.  Equivalently, the mass norm $M_K(T)$ of $T$ on each compact set $K$, is finite.  Associated with such a current $T$ is a Radon measure $\|T\|$  and a generalized tangent space $\oa T_x\in\L_pT_xX$ defined for $\|T\|$ a. a. points $x$. Recall that each $\oa T_x$ has mass norm one.  For any $p$-form $\a$ with compact support
$$
T(\a)\ =\ \int\a(\oa T)\, d\|T\|
\eqno{(\DD.1)}
$$

\Def{\DD.3}  At each point $x\in X$ let $\L(\f)$ denote the span of $G(\f)\subset \L_p TX$, 
and let 
$$
\lp \ \subset \ \L(\f)
$$
denote the convex cone  on $G(\f)$ with vertex the origin.  The $p$-vectors
$\x\in \lp$ will be called {\bf $\lp$-positive}.
\medskip

Note that $\lp$ is just the cone on $ {\rm ch}  \, G(\f)$, the convex hull of the Grassmannian.

The following Lemma is needed for a robust understanding of the definition of a 
$\f$-positive current.

\Lemma{\DD.4} {\sl The following conditions are equivalent:
\smallskip

\qquad\qquad\qquad 1)\ \ $\oa T \in \L_+(\f) \qquad\ \  \|T\| $-a.e.
\smallskip

\qquad\qquad\qquad 2)\ \ $\oa T \in {\rm ch}  \, G(\f) \qquad \|T\|$-a.e.
\smallskip

\qquad\qquad\qquad 3)\ \ $\phi(\oa T )=1 \qquad\ \ \ \   \|T\|$-a.e.
\medskip\noindent
}

The proof is provided later.

\Def{\DD.5}  A {\bf $\f$-positive  current}  is a $p$-dimensional current $T$
which is representable by integration and for which  the equivalent 
conditions of Lemma \DD.4 are satisfied.

\Prop{\DD.6} {\sl  Suppose $T$ is a compactly supported $p$-dimensional current which is representable by integration.  Then
$$
T(\f)\ \leq \ M(T)
$$
with equality if and only if $T$ is a $\f$-positive current.

Consequently, any $\f$-positive current $T_0$ with compact support is homologically mass-minimizing, i.e., 
$$
M(T_0)\ \leq \ M(T)
\eqno{(\DD.2)}
$$
for any $T=T_0+dS$ where $S$ is a $(p+1)$dimension current with compact support.  Furthermore,
equality holds in (\DD.2) if and only if $T$ is also $\f$-positive.}

\pf
Note that $T(\f)=\int\f(\oa T)d\|T\| \leq \int\|T\| = M(T)$ since $\f(\oa T)\leq \|\f\|^*\|T\|=1$.
Equality occurs if and only if $\f(\oa T)=1$  ($\|T\|$-a.e.).  This is Condition 3) in Lemma \DD.4.
The second assertion follows from the fact that $T_0(\f)=T(\f)$.  \qed

\medskip
 
 The reader may note that only Condition 3) of \DD.4 was used in this proof.
 However, it is Conditions 1) and 2) which
give a genuine understanding of $\f$-positive currents.

The fact that 
$$
M(T)\ =\ T(\f)\ =\ \int \f(\oa T)\, d\|T\|
\eqno{(\DD.3)}
$$
for all $\f$-positive currents $T$, has important implications.
\medskip

Deep results in geometric measure theory
have important applications here.

\Theorem{\DD.7} {\sl Fix a compact set $K\subset X$ and a constant $c>0$. Then the set
$\cp(\f, K,c)$ of $\f$-positive currents $T$ with $M(T)\leq c$ and $\supp(T)\subseteq K$ is compact
in the weak topology.  }
\pf    Proposition  \DD.6 easily implies that a weak limit of $\f$-positive currents is $\f$-positive.
The result then follows from standard compactness theorems for measures.\qed

\Theorem{\DD.8} {\sl Fix a compact set $K\subset X$ and a constant $c>0$. Then the set
${\cal R}(\f, K,c)$ of  rectifiable $\f$-currents $T$ with $M(T)\leq c$ and $\supp(T)\subseteq K$ is compact
in the weak topology.}
\pf
This  follows from Proposition \DD.6 and the 
Federer-Fleming weak compactness theorem for rectifiable currents [FF], [F].
\qed

\Theorem{\DD.9} {\sl  Let $T$ be a  $\f$-cycle on $X$.  Then
there is a closed subset $\Sigma\subset \supp(T)$ of Hausdorff  dimension $p-2$ such that
$M \equiv \supp(T)-\Sigma$ is a $\f$-submanifold of $X$ and 
$$
T\ =\ \sum_k n_k [M_k]
$$
where the $n_k$'s are positive integers and the $M_k$'s are the connected components of $M$.
}
\pf This is a direct consequence of Almgren's  big regularity theorem [A].\qed \medskip

We now present a dual characterization of $\f$-positive currents which will prove useful
in the next two sections.

Let $\lpp\subset \Lambda^p V$ denote the polar cone of $\lp\subset \Lambda_p V$.
By definition this is the set of $\a\in \L^pV$ such that $\a(\x)\geq0$ for all $\x\in \L_+(\f)$, or 
equivalently,
$$
\lpp = \{\a \in \L^p V : \a(\x)\geq 0 {\ \ \rm for\ all\ } \x\in
G(\f)\}.
$$
A   $p$-form $\a\in \L^pV$  is said to be {\bf  $\lpp$-positive} 
if  $\a \in \lpp$, and {\bf  strictly $\lpp$-positive} if $\a(\x)>0$      
for all $\x\in G(\f)$ (or equivalently, $\a$ belongs to the interior of
 $\lpp$).

\Remark{\DD.10} Note that $\f$ itself is strictly $\lpp$-positive,
i.e.,  an interior point of the cone $\lpp\subset \L^p V$.
 If a closed convex cone has one interior point, then there
exists a basis for the vector space consisting of interior points.
Consequently, $\L^p V$ has a basis of strictly $\lpp$-positive $p$-forms.

\medskip

If $(X,\f)$ is a calibrated manifold, the considerations and definitions
above apply to the tangent space $V=T_x X$ at each point $x\in X$. 

\Def{\DD.11}   A smooth $p$-form $\a$ on $X$ is {\bf  $\lpp$-positive (strictly
$\lpp$-positive) } if $\a$ is  $\lpp$-positive (strictly $\lpp$-positive) at each point  $x\in X$.

\Def{\DD.12} A (twisted) current $T$ of dimension $p$ is said to be {\bf  $\lp$-positive}
if  $$T(\a)  \geq  0$$ for all $\lpp$-positive $p$-forms $\a$ with compact
support.

\Theorem{\DD.13} {\sl  A current $T$ is $\lp$-positive if and only if it is
 $\f$-positive.}

\medskip

This result is proven in [HL$_3$, Prop. A.2 and Remark on page 83]. 
However, for the sake of completeness we include a proof.

\pf First assume that $T$ is representable by integration.  Then employing (\DD.1) $T$ is $\L_+(\f)$-positive if and only if 
$$
T(g\a)\ =\ \int g \a(\oa T) d\|T\| \ \geq \ 0
$$
for all functions $g\geq 0$ and all compactly supported $\lpp$-positive $p$-forms $\a$.
Equivalently, each measure $\a(\oa T)\|T\|$ is $\geq 0$ for the same set of $p$-forms.  In turn, this is equivalent to
$$
\a(\oa T)\ \geq \ 0 \qquad \|T\| {\rm -a.e.}.
$$
for all compactly supported $\lpp$-positive $p$-forms.  Finally, by the Bipolar Theorem [S]
this last condition is equivalent to the Condition 1) in the Lemma \DD.4.

It remains to prove that if $T$ is $\L_+(\f)$-positive, then $T$ is representable by integration.
For this we may assume that $T$ has compact support
in a small neighborhood $U$ of $X$, and by    Remark \DD.10, we may choose a
frame $\a_1,...,\a_N$ of smooth $p$-forms which are strictly $\lpp$-positive
on $U$.  Let $\x_1,...,\x_N$ denote the dual frame of $p$-vector fields,
i.e., $(\a_i, \x_j)\equiv\d_{ij}$ on $U$.  Every such current $T$ has a unique representation as $T=\sum_{j=1}^N u_j \x_j$ with
$u_j \in \cd'$ a distribution defined by $u_j(f) \equiv T(f\a_j)$ for all
test functions $f$. (Note that $\a=\sum_jf_j\a_j$ implies that $T(\a) =
\sum_jT(f_j\a_j) = \sum_j u_j(f_j) = (\sum_j u_j \x_j)(\sum_i f_i\a_i)
= (\sum_j u_j \x_j)(\a)$.)
Since $T$ is $\lp$-positive, each $u_j$ satisfies 
$$
u_j(f)\ \geq\ 0 \foral f\geq 0.
$$
By the Riesz Representation Theorem this proves that each $u_j$ is a
measure.  Therefore $T=\sum_j u_j \x_j$ is representable by integration.
\qed
\medskip

Now we give the proof of Lemma \DD.4.  As before $\f\in \L^pV$ is a calibration.
Let $K$ denote the unit mass ball in $\L_p V$, that is, the convex hull of 
the Grassmannian $G(p,V)\subset \L_p V$.

\Lemma {\DD.14}  $$ {\rm ch}\, G(\f) \ =\ \{\f=1\}\cap \partial K \ =\ \L_+(\f) \cap \partial K  $$

\pf  Note that:
\smallskip

a) \ \   $ {\rm ch}\, G(\f)  \ \subset\ \{\f=1\}$\ \ since $G(\f)\subset \{\f=1\}$.

\smallskip

b) \ \   ${\rm ch}  \, G(\f)  \ \subset\  K$\ \ since $G(\f)\subset G(p,V)$.

\smallskip

c)  \ \ $K\cap \{\f=1\}\ =\ \partial K\cap \{\f=1\}$ \ \ since $K\subset \{\f\leq 1\}$.

\smallskip\noindent
Hence, $ {\rm ch} \, G(\f)  \subset \{\f=1\}\cap\partial K$.

Conversely, suppose $\f(\x)=1$ and $\|\x\|=1$.  Since $\x\in K$,
$$
\x\ =\ \sum_j \lambda_j \x_j \qquad {\rm with\ each\ \ }\x_j \in G(p,V), \ {\rm each\ } \lambda>0, \ {\rm and\ }\sum_j \lambda_j=1.
$$
Hence, $1=\f(\x) = \sum \lambda_j\f(\x_j) \leq \sum\lambda_j=1$ forcing each $\f(\x_j)=1$ and therefore each $\x_j\in G(\f)$.

We have shown ${\rm ch}  \, G(\f)   \subset \partial K$, and by definition, 
${\rm ch}  \, G(\f)  \subset \L_+(\f)$.  Suppose $\x \in \partial K \cap \L_+(\f) $,  i.e., $\|\x\|=1$ and there exists   some $\lambda >0$ such that
$\lambda\x\in {\rm ch} \, G(\f)$.  We have already shown that 
 $ {\rm ch} \, G(\f)  \subset \partial K$, therefore $\|\lambda\x\|=1$, and hence 
 $\lambda =1$ proving that $\x\in{\rm ch} \, G(\f)$.
\qed

\Cor{\DD.15}  {\sl  Suppose $\x \in \L_pV$ has mass norm $\|\x\|=1$.  Then
$\x\in \L_+(\f)$ if and only if  $\f(\x)=1$  if and only if  $\x\in {\rm ch}  \, G(\f)$.}

\medskip

This is the required restatement of Lemma \DD.4

\Remark{}  Also note that the equation
$$
G(\f)\ =\ G(p,V)\cap \L_+(\f)
\eqno{(\DD.4)}
$$
follows easily from Lemma \DD.14.  This clarifies the notion of a rectifiable $\f$-current.
Namely, this proves that a rectifiable current is $\lp$-positive
if and only if it is a rectifiable $\f$-current, and eliminates a potential conflict in terminology.
\medskip

We finish this section with a lemma and corollary that are often useful.

A form $\a \in \lpp$ lies on the boundary of $\lpp$ if and only if there
exists some  $\x\in G(\f)$ with $\a(\x)=0$.  

\Lemma{\DD.16} {\sl For any  $\psi\in \L^p V$}
$$
\f-\psi \in {\rm bdy}\left\{ \lpp\right\} \ \ \ \Leftrightarrow\ \ \  \psi \leq 1 {\rm \ on\ } G(\phi) {\ \rm\ and \ }  \psi(\xi) =1 {\rm \ for\  some\   }  \xi \in G(\phi)
$$
\noindent
{\bf Proof.}  By definition $\phi-\psi \in \lpp$ if and only if $\phi(\xi)-\psi(\xi) = 1-\psi(\xi) \geq 0$ for all 
$\xi\in G(\phi)$. As remarked above $\phi-\psi$ lies in the boundary of $\lpp$ iff 
$\phi(\xi)-\psi(\xi) = 1-\psi(\xi) = 0$ at some point $\xi$.
\qed
\medskip

\Cor{\DD.17}  {\sl  For each unit vector $e\in V$, let $\f_e = e\hk(e\wedge\f) = \f\bigr|_W$ where
$W\equiv (\span e)^{\perp}$.  Then:}
$$
\f_e \in \partial \lpp \qquad {\rm if\ and\ only\ if\ }\qquad  e\in \span \x \ \ {\rm for\ some\ }\ \x\in G(\f).
$$

\pf
Note that  $\f_e = e\hk(e\wedge\f) = \f - e\wedge(e\hk \f)$ and 
$(e\wedge(e\hk \f))(\x) = |a|^2$ where $e=a+b$ with $a\in \span\x$ and $b\perp \span\x$.
Now $\f_e \in\partial \lpp $ if and only if  there exists $\x\in G(\f) $ with $|a|=1$, that is, with $e=a\in \span \x$.\qed

\Remark{}  Both $df\wedge(\nabla f\hk \f) = df\wedge d^\f f$  and  $\nabla f\hk (df \wedge \f) = \|\nabla f\|^2 \f - df \wedge d^\f f$ take values $\lpp\subset \L^pT^*X$.  Furthermore, 
\smallskip

1) \ \ $df\wedge d^\f f \ \in \  {\rm bdy}\left\{ \lpp\right\}  \qquad\Leftrightarrow\qquad \exists\ \x \in G(\f)$ tangential to the level sets of $f$.
\smallskip

2) \ \  $\nabla f \hk (df\wedge \f)\ \in \  {\rm bdy}\left\{ \lpp\right\} 
  \qquad\Leftrightarrow\qquad \exists\ \x \in G(\f)$ with $\nabla f \in \span \x$.
  \medskip
  \noindent
  Note that for some calibrations, condition 2) is true for all $f$, i.e., given a vector $n\in V$, there always exists a $p$-vector $\x \in G(\f)$ with $n\in \span \x$.

\vskip.3in


\centerline{\bf  Appendix: The reduced $\f$-Hessian.}\medskip

We assume throughout this section that $\L(\f)$ is a vector subbundle
of $\L_pTX$, and we let $\L(\f)\subset \L^p T^*X$ denote the corresponding bundle under
the metric equivalence $\L_pTX\cong \L^p T^*X$.

\Def{\DD.18}  The {\bf reduced $\f$-hessian} $\RH : C^\infty(X)\to 
\G(X, \L(\f))$ is defined to be $\ch^\f$ followed by orthogonal projection  onto the 
subbundle  $\L(\f)\subset \L^p T^*X$.  \medskip

Note that a function $f$ is $\f$-pluriharmonic if and only if $\RH (f)=0$.

Note also that if $\f$ is parallel, then $\RH = \overline d d^\f$ where $ \overline d $ denotes
the exterior derivative followed by orthogonal projection onto $\L(\f)$.

For most of the calibrations considered as examples in this paper, the image of the map
$\BM_\f:\Sym^2(TX)\to \L^p T^*X$ is contained in $\L(\f)$, or equivalently, $\RH = \ch^\f$. 
For reference, $\RH = \ch^\f$ in the following cases.
 \smallskip

(1) \ \ $\f = {1\over p!} \o^p$, the $p$th power of the K\"ahler form,

(2) \ \ $\f$ Special Lagrangian

(3) \ \ $\f$ Associative, Coassociative or Cayley

(4) \ \ $\f$ the fundamental 3-form on a simple Lie group.

\smallskip

Exceptions will be discussed at the end of this appendix.

Even when  $\RH = \ch^\f$ the following proposition is important.

\Prop{\DD.19}  {\sl  Suppose $f$ is a distribution on $X$.  Then $f$ is \fp if and only if 
$\RH(f)\equiv R$ is representable by integration and $\overrightarrow R\in \L^+(\f)$
$\|R\|$-a.e., that is, if and only if $\RH(f)$ is a $\f$-positive current.}

\medskip
The proof is similar to the proof of Theorem \DD.13 and is omitted.

\Def{\DD.20}  The $\f$-Grassmannian $G(\f)$ {\bf involves all the variables}  if, for $u\in TX$, 
the condition $u\hk \x=0$ for all $\x \in G(\f) $ implies   $u=0$
 
 \Ex{} The 2-form $\f\equiv dx_1\wedge dx_2+\lambda dx_3 \wedge dx_4$ with $|\lambda|<1$
 is a calibration on $\bbr^4$ which involves all the variables (See Section \AA), but the only
 $\x\in G(\f)$ is the $x_1,x_2$ plane so that $G(\f)$ does not involve all the variables.
 
 \Prop{\DD.21} {\sl  The operator $\RH$ is overdetermined elliptic  if and only if $G(\f)$ involves all the variables.}
 
 \pf
 We  need only consider the case $\f\in \L^p V$, where $V$ is an inner product space.
 The symbol of $\ch^\f$  at $u\in V$ is $u\wedge(u\hk \f)$.  Hence, the reduced 
 operator $\RH$ is elliptic if and only if 
 $$
( u\wedge(u\hk \f))(\x)=0\ \  \forall \ \ \x \in G(\f) \qquad\Rightarrow \qquad u=0.
 $$
 For $\x\in G(p,V)$ and $u\in V$, let $u=a+b$ with $a\in\span \x$ and $b\perp\span \x$.
 Then $(u\wedge(u\hk \f))(\x) = \f(u\wedge(u\hk \x)) = \f((a+b)\wedge( a\hk \x))
 = |a|^2\f(\x) +\f(b\wedge(a\hk\x))$.  If $\x\in G(\f)$, then
 $\f(b\wedge(a\hk\x))=0$ by the First Cousin Principle, and  $\f(\x)=1$.  Hence, 
 $(u\wedge(u\hk \f))(\x) = |a|^2=0$ if and only if $u\hk\x=0$. \qed
 
 \medskip
 
 One can easily reduce a calibration to the elliptic case.
 
 \Prop{\DD.22}  {\sl  Suppose $\f\in \L^p V$ is a calibration.  Define $W\subset V$ by 
 $$
 W^\perp\ \equiv\ \bigcap_{\x\in G(\f)}
(\span\x)^\perp 
$$
and set $\psi = \phi\bigr|_W$.  Then $\psi\in \L^pW$ is a calibration and $G(\psi)$ involves
all the variables in $W$.  Moreover, $G(\f)=G(\psi)$ and the reduced operators
$\RH$ and $\overline {\ch}^{\psi}$ agree.}
 
 \pf
 Obviously $\psi$ is a calibration and $G(\psi)\subset G(\f)$.  
 If $\x\in G(\f)$, then $\span \x\subset W$ and hence $\f(\x) = \psi(\x)$. Thus $G(\f)=G(\psi)$.
 By construction $G(\psi)$ involves all the variables in $W$.  Finally, for all $\x \in G(\f)$,
we have $\RH (f)(\x) =\tr_\x \Hess f= \overline {\ch}^{\psi}(f)(\x)$.\qed
 
 \Ex{} Let  $\Psi\in \L^4_{\bbr}\bbh^n$ be the quaternionic calibration (\AA.5) on $\bbh^n$.
 One can show that $dd^\Psi f=0$ if and only if $\Hess f=0$.
 However,
  $$
\overline{d} d^\Psi f\ =\  \overline {\ch}^{\Psi}(f)\  =\ 
 \BM_{\Psi}\left(\!\!\left({\partial^2f\over \partial q_{\a} \partial \bar{q}_{\beta}}\right)\!\!\right),
 $$
 that is, the reduced hessian is isomorphic to the quaternionic hessian
 $
 \left({\partial^2f\over \partial q_{\a} \partial \bar{q}_{\beta}}\right)
 $

\vfill\eject

\centerline{\bf \EE. Duval-Sibony Duality.}\medskip

In this section we extend the fundamental duality results established in [DS] in the complex case to calibrated manifolds $(X,\f)$. The  Duval-Sibony duality results involve plurisubharmonic functions,
pseudoconvex hulls, positive currents and Poisson-Jensen formulas.

\Def{\EE.1}  Suppose $R$ is a $(p-1)$-dimensional current on $X$.  The operator $\bdf$ is defined by
$$
(\bdf R)(f)\ =\ R(d^\f f)
$$
for all $f\in C^{\infty}_{\rm cpt}(X)$.

In other words, $\bdf : \cd_p'(X) \arr \cd_0'(X)$ is the formal adjoint of $d^\f:\ce^0(X) \arr \ce^p(X)$.
Let $\bd :\cd_{p}'(X)\arr \cd_{p-1}'(X)$ denote the {\bf boundary operator} on currents.  This is the formal adjoint of $d:\ce^{p-1} \arr \ce^p(X)$ and is related to the deRham differential on currents by
$\bd =(-1)^{n-p}d$.  The formal adjoint of $dd^\f:\ce^0(X)\arr \ce^p(X)$ is the operator
$$
\bdf \bd : \cd_p'(X) \ \arr\ \cd_0'(X)
\eqno{(\EE.1)}
$$

\Remark{}   Throughout the remainder of this section we assume that $(X,\f)$ is a non-compact connected calibrated manifold.  We also assume that $\f$ is parallel.  This assumption enables
us to use the operator $\partial_\f \partial$, but it is not necessary.  We leave it to the reader
to  verify that all of the results of this section extend to the case where $\f$ is not parallel  by replacing
 the operator  $\partial_\f \partial$ with $\ch_\f : \cd_p'(X)\to  \cd_0'(X)$, the formal adjoint of 
 $\ch^\f:\ce^0(X)\to \ce^p(X)$.  Of course, $\ch_\f$ is defined by 
 $(\ch_\f(T))(f) = T(\ch^\f(f))$ for all $f\in C^{\infty}_{\rm cpt}(X)$.

\Lemma{\EE.2. (The Support Lemma)}  {\sl  Suppose $K$ is a compact subset of $X$.  Suppose $T$ is a $\lp$-positive current with compact support in $X$.  If $\bdf \bd T$ is $\leq 0$ (a non-positive measure) on $X-\wh K$, then 
$\supp \,T \subset \wh K\cup \Core(X)$.}

\pf  Lemma \BB.2 states that for each $x\notin \wh K\cup \Core(X)$ there exists a non-negative \fp function $f$ on $X$ which is identically zero on a neighborhood of $K$  and strict at $x$.  Since $f$ is strict  at $x$, there exists a small ball $B$ about $x$ and $\e>0$ so that $dd^\f f-\e \f$ is $\lpp$-positive at each point in $B$.  By equation (\DD.3), $M(\chi_B T)= (\chi_B T, \f)$.  Therefore, 
$\e M(\chi_B T) = (\chi_B T, \e \f) \leq (\chi_BT , dd^\f f)\leq (T,dd^\f f) = (\bdf \bd T, f)\leq 0$.
This proves that $M(\chi_B T)=0$  and hence $\supp\, T\subset \wh K\cup \Core(X)$.\qed
\medskip

The case where $K=\emptyset$ is a generalization of  Proposition \BB.12.

\Cor{\EE.3} {\sl  If $T$ is a $\lp$-positive current with compact support and $\bdf\bd T\leq 0$, then}
$$
\supp\, T\ \subset \  \Core(X).
$$

When $\Core(X)=\emptyset$ we have

\Cor{\EE.4} {\sl  Suppose $(X,\f)$ is strictly convex and $K$ is $\f$-convex.  
Suppose $T$ is $\lp$-positive with compact support.  If $\supp \{\bdf\bd T\}\subset K$, then $\supp\, T\subset K$.
In particular, there are no $\lp$-positive currents which are compactly supported without boundary on  $X$.}
\medskip
.

Suppose $\ol M = M\cup \partial M$ is a compact oriented submanifold with
boundary in $X$, and that $M$ has no compact components. Let $G_x$ denote
the Green's function for $(M, \partial M)$ with singularity at $x\in M$. Let  $\mu_x$
denote harmonic measure (i.e., the Poisson kernel)  on $\partial M$ and let $[x]$ denote the 
point-mass measure at $x\in M$.  Then
$$
*_M \Delta_M G_x\ =\ \mu_x - [x]\qquad {\rm on\ \ }\ol M.
\eqno(\EE.2)
$$

If $\ol M$ is a $\f$-submanifold of a calibrated manifold $(X,\f)$, this  equation can be reformulated as a
current $\bdf \bd$-equation on $X$.

\Lemma{\EE.5}  {\sl  Suppose $M$ is a $\f$-submanifold of $X$, and that $u\in
\cd'^0(M)$ is a generalized function on $M$.  Then  }
$$
\bdf \bd(u[M])\ =\ (*_M\Delta_M u)[M].
$$

\pf
Consider the inclusion map $i:M\hookrightarrow X$.  Then, by definition,
$u[M] = i_* u$ and $(*_M\Delta_M u)[M] = i_*(*_M\Delta_M u)$.  For any test
function $f$ on $X$ we have 
$  
((\bdf \bd)(i_* u), f) \ =\ (i_*u, dd^\f f)\ =\ (u, i^*(dd^\f f))_M
$
where $(\cdot , \cdot)_M$ denotes the pairing of functions with currents on
$M$.  Proposition \AA.13 states that $i^*(dd^\f f) \ =\ *_M \Delta_M(i^* f)$
if $M$ is a $\f$-submanifold.   Finally,
$
(u,  *_M \Delta_M(i^* f))_M\ =\ ( *_M \Delta_M u, i^* f)_M
\ =\ (i_*( *_M \Delta_M u), f)_X.
$
\qed

\Cor{\EE.6} {\sl  Suppose $\ol M = M \cup \partial M$, as above, is a
$\f$-submanifold with boundary.  Then
$$
\bdf \bd (G_x[M])\ =\ \mu_x - [x]
$$
as a current equation on $X$.}
\medskip

Assume $K$ is a compact subset of $X$, and let $\cm_K$ denote the set of
probability measures with support in $K$.  

\Def{\EE.7}  If $T_x$ is a $\lp$-positive current with compact support and
$T_x$ satisfies:
$$
\bdf \bd \,  T_x\ =\ \mu_x -[x]
\eqno(\EE.3)
$$
with $\mu_x\in \cm_K$, then: $T_x$ is a {\bf Green's current  for }
$(K,x)$, $\mu_x$ is a {\bf Poisson-Jensen measure} for $(K,x)$, and the
equation (\EE.3) is the {\bf Poisson-Jensen equation}.

\medskip

\Theorem{\EE.8} {\sl  Suppose $X$ is strictly $\f$-convex, $K$ is a
compact subset of $X$, and $x\in X-K$. Then there exists a Green's current
$T_x$ for $(K,x)$ if and only if $x\in \PH K$.
}

\medskip
To prove this we begin with the following.

\Prop{\EE.9}  {\sl  Suppose $(X,\f)$ is non-compact   calibrated manifold. 
If there exists a Green's current for $(K,x)$, then $x\in \PH K$.
}
\medskip

\pf This follows immediately from Lemma \EE.2 since $x\in\supp\, T_x$.\qed

\medskip
\noindent
{\bf Second Proof.}
Since $\bdf \bd \, T_x \ =\ \mu_x-[x]$, we have $\int f \mu_x - f(x) \ =\
(T_x, d d^\f f)$ for all $f\in C^\infty(X)$.
If $f$ is \fp, this implies that $f(x)\ \leq\  \int f \mu_x\ \leq \ \sup_K
f$, since $\mu_x\in \cm_K$. Thus $x\in \PH K$.\qed
\medskip

The set $\cp_X \equiv \fpsh\subset C^\infty(X)$  of all
\fp-functions on $X$ is clearly a closed convex cone in
$C^\infty(X)$.  Let 
$$
C_X\ \equiv\ \{u\in \cd'_{0, {\rm cpt}}(X) : u=\bdf \bd \, T \ {\rm for\ some\ } \lp{\rm-positive\ }T\in
\cd'_{p, {\rm cpt}}(X)  \}. 
\eqno{(\EE.4)}
$$
This is a convex cone in $ \cd'_{0, {\rm cpt}}(X)$.

\Lemma{\EE.10} {\sl  Suppose $X$ is non-compact with calibration $\f$.
Then $\cp$ is the polar of $C$,
that is,
$$
\cp \ =\ C^0\  \equiv \ \{f\in C^\infty(X) \ : \  (u,f)\geq 0 \ \  \forall u
\in C\}. $$
}
\pf
Consider $u=\bdf \bd(\delta_x \x)$, with $\x\in G_x(\f)$.  Clearly $u\in C$. 
If $f\in C^\infty(X)$ belongs to $C^0$, then 
$0\leq (u,f) = (\bdf \bd (\delta_x \x), f) \ =\ (\delta_x \x, d d^\f f)\ =\
(dd^\f f)_x(\x).$  Hence $C^0\subseteq \cp$.

Conversely, if $f\in\cp$, then for all $u\in C$, $(u,f) \ =\ (\bdf \bd   T, f)\ =\ (T, dd^\f
f)\geq 0$, since $T$ is $\lp$-positive. This proves that $\cp\subseteq C^0$.
\qed

\Lemma{\EE.11}  {\sl  If $X$ is strictly $\f$-convex, then the convex cone
$C\subset \cd'_{0, {\rm cpt}}(X)$ is closed.} 

\pf
It suffices to show that $C\cap \cd'_{0, K}(X)$ is closed for an exhaustive
family of compact subsets $K\subset X$.  We may assume $K$ is $\f$-convex.
Suppose $u_j$ converges to $u$ in $C\cap \cd'_{0, K}(X)$ with each $u_j\in C$,
i.e., $u_j = \bdf \bd T_j$ where $T_j$ is a $\lp$-positive current with
compact support.  By Corollary \EE.4  the support of each $T_j$ is
contained in $K$.  Consider a strictly \fp function $f$ on $X$.  Pick
$\epsilon>0$ so that $dd^\f f -\epsilon \f$ is $\lpp$-positive at each
point of $K$.  Then 
$
\epsilon M(T_j) \ =\ (T_j, \epsilon \f) \ \leq\ (T_j, dd^\f f)\ =\ (\bdf \bd 
T_j, f)\ =\ (u_j,f) 
$
which converges to $(u,f)$.  Therefore the masses $M(T_j)$ are bounded.
By compactness there exists a weakly convergent subsequence $T_j\to T$.
Now $\supp\, T\subset K$ and $T$ must be $\lp$-positive. 
Hence $u=\bdf \bd T\in C\cap \cd'_{0, K}(X)$.  
This proves that $C\cap \cd'_{0, K}(X)$ is closed for each compact set $K$
which is $\f$-convex. \qed

\Cor{\EE.12}  {\sl  Suppose $X$ is strictly $\f$-convex.  Then
$$
C \ =\ \cp^0.
$$
Equivalently, the equation 
$$
\bdf \bd  T\ =\ u
$$
 has a solution $T$ which is a
$\lp$-positive current with compact support if and only if}
$$
0\ \leq \  u(f)\qquad {\sl for\ all\ \ \ } f\in\fpsh
$$
\pf 
Apply the Bipolar Theorem.\qed

\medskip
\noindent
{\bf Proof of Theorem \EE.8.}
Suppose there does not exist a Green's current for $(K,x)$, that is,
suppose $\cm_K-[x]$ is disjoint from the cone $C$.  By the Hahn-Banach
Theorem (note that $\cm_K-[x]$ is a compact convex set) there exist $f\in
C^0=\cp$ with $f$, considered as a linear functional on $\cd'_{0, {\rm cpt}}(X)$, satisfying $u(f) \leq -\e <0$ for all $u\in (\cm_K-[x])$.   That is, $\int f d\,\mu - f(x) \leq -\epsilon <0$ for all $\mu\in
\cm_K$. Consequently, 
$$
\sup_K f \ =\ \sup_{\mu\in \cm_K} \int f  d\,\mu\ \leq\ f(x) -\epsilon
$$
or $\sup_K f +\epsilon \leq f(x)$ so that $x\notin \PH K$.\qed
\medskip

One could define the ``Poisson-Jensen hull'' of a compact set $K$ to be the set of points $x$
for which there exists a Poisson-Jensen measure $\mu_x$ and a Green's current $T_x$ satisfying
(\EE.3).
Then Proposition \EE.9 states that on any (non-compact) calibrated manifold $(X,\f)$, the 
Poisson-Jensen hull of a compact set is contained in the \fp hull, while Theorem \EE.8
states that the two hulls are equal if $(X,\f)$ is strictly convex.

The next ``hull"  obviously contains the Poisson-Jensen hull.

\Def{\EE.13}  The {\bf current hull} of a compact subset $K\subset X$ is the union
$$
\wt K\ \equiv\ \bigcup_{T\in \cp(K)}
\supp\, T
$$
where $\cp(K)$ consists of all $\lp$-positive currents with compact support on $X$ 
satisfying $\bdf \bd  T\leq 0$ on $X-K$.

\Lemma{\EE.14}  {\sl  If $(X,\f)$ is strictly $\f$-convex, then $\wt K = \wh K$.}

\pf
The support Lemma \EE.2 states that 
$$
\wt K\ \subset \ \wh K \cup \Core(X).
$$
for any calibrated manifold.  Now $\wt K$ contains the Poisson-Jensen hull
which equals $\wh K$ if $(X,\f)$ is strictly \fc  by Theorem 5.8.\qed

Suppose $(X,\f)$ is non-compact and connected.

\Def{\EE.15}
  An open subset $\Omega \subset X$ is   {\bf \fc  relative to $X$} if $K\subset\subset \Omega$
  implies ${\wh K}_X\subset\subset \Omega$.\medskip
  
  Note that if $X$  is \fc this condition implies that $\Omega$ is \fc since
 ${\wh K}_\O\subseteq{\wh K}_X$.
 Moreover, if $X$ is strictly \fc, then $\Core(\O)\subseteq\Core(X)$ is empty so that $\O$ is
 strictly \fc (by Proposition 2.13).

\Prop{\EE.16}  {\sl  Suppose $(X,\f)$ is strictly \fc.  An open subset  $\O^{\rm open}\subset X$ is \fc relative to $X$ if and only if  ${\cal PSH}(X,\f)$ is dense in ${\cal PSH}(\O,\f)$.}

\pf Let $L:C^\infty(X)\to C^\infty(\O)$ denote restriction. The adjoint $L^*:\cd_{0,{\rm cpt}}'(\O) \to 
\cd_{0,{\rm cpt}}'(\O)$ is inclusion.  Suppose $v\in (L^*)^{-1}(C_X)$. i.e., $v\in \cd_{0,{\rm cpt}}'(\O)$ with $v=\partial_\f\partial T$ for some $\lp$-positive current $T$ compactly supported in $X$.  Then $K\equiv \supp\, v$ satisfies $\wt K_X =\wh K_X$ by Lemma \EE.14. 
 Hence, $\wh K_X\subset \O$ implies    $\supp\, v \subset \O$ or that $v\in C_\O$.  
 This proves that $\O$ is \fc relative to $X$ if and only if 
$$
 (L^*)^{-1}(C_X) \ =\  C_\O.
\eqno{(\EE.5)}
$$
By Corollary \EE.12 we may replace $C_X$ by $\cp_X^0 = C_X$.  In general, 
$[L(\cp_X)]^0 = (L^*)^{-1}(\cp_X^0)$, so that (\EE.5) is equivalent to
$$
\left[ L(\cp_X) \right]^0 \ =\ C_\O.
\eqno{(\EE.6)}
$$
By Lemma \EE.10, $\cp_\O = C_\O^0$.  Hence (\EE.6) is equivalent to 
$$
\overline{L(\cp_X)} \ =\ \cp_\O
\eqno{(\EE.7)}
$$
\qed

\vskip .3in



\vskip .3in

\centerline{\bf \FF.  \   $\f$-Free Submanifolds}

\medskip

Suppose $(X,\f)$ is a calibrated manifold.  A plane $\x$ is said to be {\sl tangential} to a submanifold $M$ if $\span \x \subset T_xM$.

\Def{\FF.1}  A closed submanifold $M\subset X$ is {\sl  \totr} if there are no $\f$-planes $\x\in G(\f)$ which are tangential to $M$.
If the restriction of the calibration $\f$ to $M$ vanishes,   $M$ is called  {\sl  $\f$-isotropic}.

\medskip

Note that  $\f$-isotropic submanifolds are  \totr.
Each submanifold of dimension strictly less than the degree of $\f$ is $\f$-isotropic and hence
automatically \totr.  Furthermore, in dimension $p$ the generic local submanifold is \totr.
Depending on the geometry, this may continue through a range of dimensions greater than $p$.

\Theorem {\FF.2}  {\sl Suppose $M$ is a closed submanifold of $(X,\f)$ and let 
$\dis(x) \equiv {1\over2}{\rm dist}(x,M)^2$ denote half the square of the distance to $M$.  Then
$M$ is \totr \ if and only if 
the function $\dis$ is strictly \fp at each point in $M$ (and hence in a neighborhood of $M$).
}

\pf  We begin with the following.

\Lemma{\FF.3}  {\sl  Fix $x\in M$ and let $P_N:T_xX\to N$ denote orthogonal projection onto the normal plane   of $M$ at $x$.  Then for each $\x \in G(\f)$ one has}
$$
\{\BM_\f(\Hess_x \dis)\}(\x) \ =\ \langle P_N, P_\x\rangle
\eqno {(\FF.1)}
$$
\pf  By Lemma 1.15b)  
$$
\{\BM_\f(\Hess_x f)\}(\x) \ =\ \langle \Hess_x f, P_\x\rangle
\eqno {(\FF.2)}
$$
for any function $f$.
The lemma then follows from the assertion that 
$$
\Hess_x \dis \ =\ P_N.
\eqno{(\FF.3)}
$$
To see this we first note that the Hessian of any function $f$  can be written
$$
\Hess f (V,W) \ =\ \langle V, \nabla_W( \nabla f)\rangle 
\eqno{(\FF.4)}
$$
for all  $V,W\in T_xX$.
It follows that if $\nabla f =0$ on the submanifold $M$, then $T_xM \subset {\rm Null}(\Hess_x f)$.
Thus, with respect to the decomposition $T_xX=T_xM \oplus N$ we have
$$
\Hess_x \dis \ =\ \left(\matrix { 0&0\cr 0&A\cr }\right)
$$
and it remains to show that $A$ is the identity. To see this, set $\delta(x)= \dist(x, M)$  
and note that $\nabla \delta = n$ is a smooth unit-length vector field near (but not on)
$M$ whose integral curves are geodesics emanating from $M$.   Hence, $$\nabla_n (\nabla \dis) =
\nabla_n (\nabla {1\over2}\delta^2) = \nabla_n (\delta n) = n + \delta \nabla_n n= n.$$
Taking limits along normal geodesics  down to $M$ gives the result.
\qed
\medskip

Theorem \FF.2 now follows from the fact that 
$$
\langle P_N, P_\x\rangle\ \geq\ 0 \ {\rm with \ equality \ iff }\ \span\x \subset N^\perp = T_xM.
$$
\qed

The existence of \totr \  submanifolds insures the existence of lots of strictly \fc domains in $(X,\f)$.

\Theorem{\FF.4}  {\sl  Suppose $M$ is a \totr\ submanifold of $(X,\f)$.
  Then there exists a fundamental neighborhood system $\cf(M)$ of $M$ such that:
  
  \smallskip
  
  (a)\  $M$ is a deformation retract of each $U\in \cf(M)$.

  \smallskip
  
  (b)\  Each neighborhood $U\in \cf(M)$ is strictly \fc.

  \smallskip
  
  (c)\  \ ${\cal PSH}(V,\f)$ is dense in  ${\cal PSH}(V,\f)$ if $U\subset V$ and $V,U\in\cf(M)$.

  \smallskip
  
  (d)\  Each  compact set $K\subset M$ is ${\cal PSH}(U,\f)$-convex for each $U\in \cf(M)$.
  }

\pf   We construct tubular neighborhoods of $M$ as follows.  Let  $\epsilon \in C^\infty(M)$ 
be a smooth function which vanishes at infinity and has the property that  for each $x
\in M$ the ball $\{y\in X:{1\over2}\dist(y,x)^2\leq \epsilon(x)\}$ is compact and geodesically convex.
Assume also that $\e$ is sufficiently small so that the exponential map gives a diffeomorphism
$$
\exp : N_\e \ \arr\ U_\e
$$
from the open set $N_\e$  in the normal bundle $N$ defined by
${1\over2} \|n_x\|^2<\e(x)$ to the  neighborhood
$$
U_\e \ =\ \{x\in X : \dis(x) < \e(x)\}.
\eqno{(\FF.5)}
$$
of $M$ in $X$.
Each $U_\e$ admits a deformation retraction onto $M$.

By Theorem \FF.2 the function $\dis={1\over2}{\rm dist}(\cdot, M)^2$ is strictly \fp\ on a neighborhood of $M$, which we can assume to be $W$.  We impose the following additional condition on the function $\e\in C^\infty(W)$.
$$
\dis - t\e \ \ \  {\rm is\ strictly \ } \f-{\rm plurisubharmonic \ on\ } \ W \ \ {\rm for }\ \ 0\leq t \leq 1.
\eqno{(\FF.6)}
$$
Since (\FF.6) is valid as long as $\e$ and its first and second derivatives vanish sufficiently fast
at infinity, it is easy to see that the family $\cf(M)$ of neighborhoods $U_\e$ constructed above with $\e$ satisfying (\FF.6) is a fundamental neighborhood system for $M$.

Obviously, the function $\psi \equiv (\e-\dis)^{-1}$ is a proper exhaustion for $U_\e$.  To prove (b), recall that if $g$ is a positive concave function, then $1/g$ is convex, or more directly, calculate that
$$
\Hess \psi \ =\ \psi^2\Hess(\dis-\e) +\psi^3\nabla(\e-\dis)\circ \nabla(\e-\dis).
\eqno{(\FF.7)}
$$
Applying $\BM_\f$ to (\FF.7) proves that $\psi$ is strictly \fp\ on $\{\psi>0\}=U_\e$.

To prove parts (c) and (d) one uses Proposition \EE.16 and argues exactly as on page 302 of 
[HW$_1$].       \qed

\medskip

\Ex{\FF.5}  As mentioned above, Theorem \FF.4 exhibits a rich family of \fc \ domains in $(X,\f)$.
For example, let $M\subset X$ be {\sl any submanifold of dimension} $< p = \deg \f$. Then by
\FF.4, $M$ has a fundamental system of neighborhoods each of which is a strictly $\f$-convex domain homotopy equivalent to $M$.

\Ex{\FF.6} Interesting examples occur in all the calibrated geometries examined in depth in [HL$_3$].
Suppose for instance that $X$ is a Calabi-Yau manifold with special Lagrangian calibration
$\f$.  Then {\sl any complex submanifold $Y\subset X$ is \totr.}  It follows that any smooth submanifold
of $Y$ is also \totr.
\medskip

We now consider  the following two classes of subsets of $(X,\f)$.
\smallskip

\qquad (1)\ Closed subsets $A$ of \totr submanifolds.
\smallskip

\qquad (2)\ Zero sets  of non-negative strictly \fp\ functions $f$.
\smallskip

These two classes are basically the same, as described in the following two propositions.

\Prop{\FF.7}  {\sl  Suppose $A$ is a closed subset of a \totr submanifold $M$ of $X$.  Then
there exists a non-negative function $f\in C^\infty(X)$ with
\smallskip

(a) \ \ $A=\{x\in X: f(x)=0\}$
\smallskip

(b) \ \ $f$ is strictly \fp\ at each point in $M$ (and hence in a neighborhood 

\qquad  of $M$ in $X$).
}

\pf  Since $M$ is a closed submanifold, the function  $\dis$ in Theorem \FF.2 can be extended to
$h\in C^\infty(X)$ which agrees with $\dis$ in a neighborhood of $M$ and satisfies
$$
h\geq 0 \and \{h=0\}= M.
$$

Choose $\psi\in C^\infty(X)$ with $\psi\geq0$ and $A=\{x\in X: \psi=0\}$.  Now choose $\e\in
 C^\infty(X)$  with $\e(x)>0$ for all $x\in X$, and with $\e$ and its derivatives sufficiently small so that $f\equiv h+\e\psi$ is strictly \fp \ on $M$.
\qed

\Prop{\FF.8} {\sl  Suppose $f  \in C^\infty(X)$ is a non-negative function which is strictly 
\fp at each point in $A \equiv \{x\in X:f(x)=0\}$. Given a point $x\in A$ there exists a neighborhood $U$ of $x$ and a proper \totr \ submanifold $M$ of $U$ such that $A\cap U\subset M$.
}
\pf
Given $x\in A$ we may choose geodesic normal  coordinates $(z,y)$ in a neighborhood $U$ at $x$  so that
$$
\Hess_x f\ =\ \left( \matrix{0&0\cr0&\Lambda\cr}\right)\eqno{(\FF.8)}
$$
where $\Lambda$ is the diagonal matrix diag$\{\lambda_1,...,\lambda_r\}$, $r$ is the rank
of $\Hess_x f$, and $\lambda_j\neq 0$ for $j=1,...,r$.
Set 
$$
M\ =\  \left \{w\in U: {\partial f\over\partial y_1}=\cdots ={\partial f\over\partial y_r}=0\right\}.
$$
Since $\nabla{\partial f\over\partial y_1},..., \nabla{\partial f\over\partial y_r}$ are linearly independent at $x$, $M$ is a codimension $r$ submanifold locally near $x$.  

Note that ker$(\Hess_x f)=T_xM$. It remains to show that ker$(\Hess_x f)$ is totally real,
 since  if $M$ is \totr at $x$, then $M$ is \totr in a neighborhood of $x$. This  is proved in Lemma \FF.9 below.\qed\medskip
 
 \Lemma{\FF.9}  {\sl  Suppose $f$ is  strictly \fp at $x\in X$.   Then $\ker(\Hess_x f)\subset T_xX$
 is \totr.     }

\pf
If  $\ker(\Hess_x f)\subset T_xX$ is not \totr, there exists $\x\in G(\f)$ with 
$(\Hess_x  f)\bigr|_{\span \x}=0$.  Consequently, 
$d d^\f f(\x) = \BM_\f(\Hess_x  f)(\x) = \tr_\x( \Hess_x f)=0$, and $f$ is not strict at $x$.
\qed

\Remark{\FF.10}  Parts (b), (c) and (d) of Theorem \FF.4 can be generalized as follows.
Suppose $M=\{f=0\}$ is the zero set of a non-negative strictly \fp function $f$ on $(X,\f)$.
Then there exists a fundamental neighborhood system $\cf(M)$ of $M$ satisfying (b), (c) and (d) of Theorem \FF.4.  The neighborhoods $U_\e\in\cf(M)$ are defined by 
$U_\e=\{x\in X: f(x)<\e(x)\}$ where $\e>0$ is a $C^\infty$ function on $X$ vanishing at infinity along with  its first and second derivatives so that $f-\e$ remains strictly \fp.  The proofs of (b), (c) and (d)
are essentially the same as in Theorem \FF.4.

\medskip

We conclude with a the following general observation.

\Prop{\FF.11} {\sl Let $M$ be a submanifold of $(X,\f)$ and $f$ a smooth function defined on
a neighborhood of $M$ such that:
\smallskip

(1) $\nabla f \equiv 0$ on $M$, and \smallskip

(2)  $f$ is strictly \fp at all points of $M$. \smallskip

\noindent
Then $M$ is \totr.}

\pf  By (\FF.4) we see that $TM\subseteq \ker(\Hess f)$ at all points of $M$.  We then apply Lemma \FF.9.\qed\medskip

\Cor{\FF.12} {\sl Let $f$ be a non-negative, real analytic function on $(X,\f)$ and consider the real analytic subvariety $Z\equiv \{f=0\}$.  If $f$ is strictly \fp at points of $Z$, then each stratum of $Z$
is \totr.}

\vfill\eject



\centerline{\bf \GG.  \   Hodge Manifolds}

\medskip

\def\IH{\wt{H}}

In this section we pose some highly speculative questions for calibrated manifolds
in the spirit of those posed  in the complex case (cf. [HK, p.58], [L$_{4,5}$]).
Assume that $(X,\f)$ is a compact calibrated $n$-manifold with a 
parallel calibration $\f$ of degree $p$.  Let
$$
\psi\ =\ *\f
\eqno{(\GG.1)}
$$
denote the dual calibration.  Note that a $\f$-submanifold or, more generally, any $\f$-cycle on $X$
is a current of dimension $p$ and degree $n-p$.  By contrast a $\psi$-submanifold or $\psi$-cycle is a current of dimension $n-p$ and degree $p$.  Denote by ${\IH}^p(X,\bbz)$  the image of  the map
${H}^p(X,\bbz) \to {H}^p(X,\bbr)$.
\Def{\GG.1} If the de Rham class of the calibration $\f$ lies in ${\IH}^p(X,\bbz)$, i.e., if $\f$ has integral periods, then $(X,\f)$ will be referred to as a {\bf Hodge manifold}.

\Remark{}  If $(X,\o)$ is a K\"ahler manifold, then this coincides with standard
terminology.  The Kodaira Embedding Theorem states that in this case each Hodge
manifold is projective algebraic with $N\o-[H]=d\a$, where $H$ is a hyperplane section,
$N$ a positive integer, and $\a$  a current of degree 1.

\medskip\noindent
\HoQu{(for the class of $\f$)} 
Suppose $(X,\f)$ is a Hodge manifold.  When does there exist a $*\f$-cycle $T$ cohomologous
to $N\f$ for some positive integer $N$, i.e.,
$$
N\f-T\ =\ d\a
\eqno{(\GG.2)}
$$
for some current $\a$ of degree $p-1$?
\medskip

Recall that a  $*\f$-cycle is automatically $*\f$-positive, so this is, more precisely,  the ``Hodge Question with Positivity'' for $\f$.

\Remark{}  If   equation (\GG.2) (called the spark equation) has a solution, then $\a$ determines a differential character on $X$.  (See [HLZ] for more details.)

\Ex{\GG.3}  In [L$_3$] an  example is constructed of a Hodge manifold $(X,\f)$ for which no such cycle exists.
More specifically, a parallel self-dual 4-form $\f$ of comass 1 is constructed  on a flat torus $X$ of dimension 8 with the property that  $[\f]\in H^4(X,\bbr)$ is an integral class, but there exist no ${\f}$-cycles  whatsoever on $X$.

\Ex{\GG.4} Consider  the fundamental bi-invariant 3-form $\O$ on a compact simple Lie group $G$,
 normalized to be the generator of $H^3(G,\bbz)\cong \bbz$. Then $(G,\O)$ is a Hodge manfiold, and R. Bryant [B] has shown  that, indeed, $\O$ is always cohomologous to a $*\O$-cycle.  These $*\O$-cycles are always sums of  singular semi-analytic subvarieties congruent to irreducible components of the cut-locus of the exponential map.
\medskip

For a general class in ${\IH}^p(X,\bbz)$ to be represented  by a $*\f$-cycle (or for a class in 
${\IH}_p(X,\bbz)\cong {\IH}^{n-p}(X,\bbz)$ to be represented  by a $\f$-cycle), there is a natural
necessary condition coming from the Hodge decomposition. Note that 
since $\f$ and $*\f$ are parallel, the subspaces $\L_x(\f)\equiv \span G(\f) \subset \L_pT_xX$
and $\L_x(*\f)\equiv \span G(*\f) \subset \L_{n-p}T_xX$, under metric equivalence $TX\cong T^*X$, 
define  {\sl parallel subbundles}
$$
\L(\f)\ \subset \L^pT^*X   \and  \L(*\f)\ \subset \L^{n-p}T^*X 
$$
The orthogonal projections $P_{\L(\f)}$ and $P_{\L(*\f)}$ onto these subbundles are {\sl parallel operators on forms}.
It was proved by Chern [Ch] that any such operator commutes with harmonic projection.
Recall that the Hodge decomposition: 
$\ce^p(X) = \bbh^p(X) \oplus {\rm Image }(d)\oplus {\rm Image }(d^*)$, is a $C^\infty$-decomposition, and therefore induces a corresponding decomposition of currents:
${\ce'}_p(X) = \bbh_p(X) \oplus {\rm Image }(\partial)\oplus {\rm Image }(\partial^*)$.

\Def{\GG.5}  A $p$-dimensional current $T$, representable by integration,  is said to be of 
{\bf type} $\Lambda(\f)$ if 
$\overrightarrow T_x \in \L(\f) \subset \L^pT_xX$ for $\| T \|$-a.a. $x$.  \medskip

This definition   extends to arbitrary currents $T$of dimension $p$ by requiring that $T(\psi)=0$ for
all smooth $p$-forms $\psi$ such that $\psi\bigr|_{\L(\f)}=0$.

\Prop{\GG.6}  {\sl  If a class $c\in \IH_p(X,\bbz)$ is represented by a current of type $\L(\f)$, then the harmonic representative of $c$ must be of type $\L(\f)$.
}

\pf
Any $\L(\f)$-current is fixed by the parallel bundle projection, and that projection commutes
with the harmonic projector. \qed

\Def{\GG.7}  A {\bf  $\L(\f)$-cycle} is a $d$-closed, $p$-dimensional locally rectifiable current of 
type $\L(\f)$.

\medskip\noindent
I.  \HoQu{}  Suppose  $u     \in  \IH_p(X,\bbz)$ is a class whose harmonic representative
is of type $\L(\f)$.  When does there exist an integer $N$ and a $\L(\f)$-cycle $T$ with
$T\in Nc$?

\Remark{\GG.8} Example \GG.3 above gives a parallel 
calibration $\f$ on a flat 8-dimensional torus $X$and an integral  class $c\in \wt {H}_4(X, \bbz)$
of type $\L(\f)$ for which no such current exists. 

\Remark{\GG.9} The Hodge Question  is  a direct generalization of the
standard Hodge Conjecture for algebraic cycles on a complex projective manifold,
since we know from [HS],  [Sh] and [Alex]  that for $\f = \o^p/p!$ ($\o$ =  the K\"ahler form), any 
$\L(\f)$-cycle is an algebraic $p$-cycle.
\medskip

Any locally finite integer sum of $\f$-cycles is a $\L(\f)$-cycle.  However, the converse is
completely open outside of the K\"ahler case.  Moreover, 
even though it holds in the K\"ahler case (cf. Remark \GG.9),  there is  no proof of this fact by the standard methods of regularity
 in Geometric Measure Theory.

Before trying to prove that a general $\L(\f)$-cycle is a sum of $\f$-cycles, one
would like the calibration $\f$ to have the algebraic property displayed in the next remark.

\Remark{\GG.10}   Equation (\DD.4) says that 
$
G(p, T_xX)\cap \L_+(\f)\ =\ G(\f)
$
so that $\f$-cycles and $\lp$-cycles are the same thing.  Most parallel calibrations(see [HL$_3$, p. 68]) 
are known to satisfy
$$
G(p, T_xX)\cap \L(\f)\ =\ G(\f) \cup (-G(\f)).
$$
In this case $T$ is a $\lp$-cycle if and only if $\pm \overrightarrow {T_x} \in G_x(\f)$ for
$\|T\|$-a. a. $x$.
Consequently, $T$ decomposes inito $T^+-T^-$ with both $\overrightarrow {T_x}^{\pm}  \in G_x(\f)$,
but, even in the K\"ahler case, one can not show directly that $T^+$ and $T^-$ are $d$-closed. 

\medskip

There are versions of the Hodge Question involving ``positivity'' which may have more hope.
For example:

\medskip\noindent
II. \HoQu{(with positivity)}
Suppose $c\in  \IH_p(X,\bbz)$ is a class whose harmonic representative is strictly $\lp$-positive.
When does there exist an integer $N$ and a $\f$-cycle $T$ with $T\in Nc$?

\Remark{\GG.11}  If the current $*\f$ (of dimension $p$)  is strictly $\lp$-positive, then for any form 
$\psi$ of type $\L(\f)$, there exists an integer $\ell$ such that $\psi+\ell(*\phi)$ is strictly
$\lp$-positive.  This applies for example to the harmonic representative of $c$ in Hodge Question II.
Consequently, one can see that {\sl if $(X, *\f)$ is a Hodge manifold with a solution to 
(\GG.2), then the Hodge Questions I and II are equivalent.}

\Remark{\GG.12}  The point of  Hodge Question II  is that one is asking for a $\f$-cycle
which is  automatically $\lp$-positive and therefore satisfies the strong regularity theorem \DD.9.

\vskip .3in



\vskip .3in

\centerline{\bf \HH.  \  Boundaries of $\f$-submanifolds.  }

\medskip

In this section we take up the following general question.  Suppose $(X,\f)$ is a non-compact strictly \fc manifold.  Given a compact oriented submanifold $\G\subset X$ of dimension $p-1$, when does there exist a $\f$-submanifold $M$ with boundary $\G$?  More generally, when does there exist
a  $\lp$-positive current $T$ with $\partial T=\G$?

\Theorem{\HH.1} {\sl Suppose $\f$ is exact.  Given $S\in \ce'_{p-1}(X)$, there  exists a $\lp$-positive current $T\in \ce'_p(X)$ with $S=\partial T$ if and only if }
$$
\int_S\a \ \geq\ 0\qquad {\rm for\ all\ }  \a \in \ce^{p-1}(X) \ {\rm such \  that \ } d\a \ {\rm is\ } \ \lpp{\rm -positive}
$$

\pf
Consider the following convex cones.
$$
A\ =\ \{\a \in \ce^{p-1}(X) : \  d\a \ {\rm is\ } \ \lpp{\rm -positive}\}
$$
$$
B\ =\ \{S\in \ce'_{p-1}(X) : \ S=\partial T \ {\rm for\ some\ } \lp{\rm -positive} \ T\in \ce'_p(X)\}
$$

If $\a\in A$ and $S\in B$, then
$$
S(\a)=\partial T(\a) = T(d\a)\geq 0,
$$
that is,
$$
A\ \subset \ B^0 \and B\ \subset \ A^0
$$
where $B^0$ denotes the polar of $B$.  If $\x\in G_x(\f)$, then $T=\delta_x \x$ is $\lp$-positive, so that $\partial(\delta_x \x)\in B$. Therefore, if $\a\in B^0$, then $0\leq \partial(\delta_x \x)(\a)
=(\delta_x \x)(d\a) = (d\a)_x(\x)$.  This proves that $B^0\subset A$, and hence $A=B^0$.
(In particular, note that $A$ is closed.)  Theorem \HH.1 is just the statement that
$B=A^0$.  Now since $A=B^0$, the Bipolar Theorem states that $\overline B = A^0$.
Thus it remains to show that $B$ is closed.

Suppose $S_j\in B$ and $S_j\to S$ in $\ce'_{p-1}(X)$.  Then $S_j=\partial T_j$ for some 
$T_j$ which is $\lp$-positive.  The calibration $\f$ is exact, i.e., $\f=d\eta$ for some
 $\eta\in\ce^{p-1}(X)$.  Therefore
 $$
 M(T_j)\ =\ T_j(\f)\ =\ T_j(d\eta)\ =\ (\partial T_j)(\eta)\ = \ S_j(\eta)\ \arr \ S(\eta).
 $$
 In particular, there exists a constant $C$ such that $M(T_j)\leq C$ for all $j$.
 By Lemma \EE.2
 $$
 \supp T_j\ \subset\ \wh{\supp S_j}
 $$
for each $j$.  Pick a compact subset $K$ with $\supp S_j \subset K$ for all $j$.
Then
$$
 \supp T_j\ \subset\ \wh{K}\qquad {\rm for\ all\ } j.
$$

This proves that $\{T_j\}$ is a precompact set in $\ce'_p(X)$.  Therefore, there exists a convergent
subsequence $T_j\to T$ in $\ce'_p(X)$. Obviously, $\partial T=S$ and $T$ is $\lp$-positive.
Hence, $S\in B$.
\qed

\Remark{\HH.2}  The same proof combined with the Federer-Fleming compactness theorem for integral currents proves the following.  Let $\ccr_p(X)$ denote the compactly supported rectifiable currents of dimension $p$ on $X$. Then, if $\f$ is exact, the set
$$
B_{\rm rect} \
 \equiv\ \{ \G\in \ccr_{p-1}(X) : \ S=\partial T \ {\rm for\ some\ } \lp{\rm -positive} \ T\in \ccr_p(X)\}\}
$$
is weakly closed in $\ccr_{p-1}(X)$.


\vskip .3in

\centerline{\bf \II.  \  $\f$-Flat Hypersurfaces and Functions which are $\f$-Pluriharmonic mod $d$. }

\medskip

The $\f$-pluriharmonic functions are the closest thing to holomorphic functions on a calibrated manifold $(X,\f)$.  Usually there are very few $\f$-pluriharmonic functions. An attempt has been made in this paper to remedy this situation by emphasizing the \fp functions. By comparison these functions exist in abundance.  For some purposes another extension of the concept of 
$\f$-pluriharmonic functions is more useful --- namely the {\sl $\f$-pluriharmonic functions mod $d$}.

This section is, for the most part, a straightforward extension of the results of Lei Fu [Fu]
from the special Lagrangian case to the general calibrated manifold $(X,\f)$.

\Def{\II.1}  A function $f\in C^\infty(X)$ is {\bf \fh mod $d$} if 
$$
dd^\f f \ =\ df\wedge\a_f + \s_f
\eqno{(\II.1)}
$$
for some 1-form $\a_f$ and some $p$-form $\s_f$ of type $\L(\f)^\perp$, i.e.,
$\s_f(\x)=0$ for all $\x\in G(\f)$.
\medskip

If $f$ is \fh mod $d$, then $\lambda f$, $\lambda\in\bbr$, is also \fh mod $d$.  However, the sum of two such functions need not be \fh mod $d$.

\Prop{\II.2}  {\sl   Suppose that $df$ never vanishes so that $\ch \equiv \ker df$ defines a hypersurface foliation.  The condition that $f$ be \fh mod $d$ is independent of the function defining 
the foliation $\ch$.}

\pf
Recall that locally $f$ and $g$ define the same foliation $\ch$ if and only if $g = \chi(f)$ for some
function $\chi:\bbr\to\bbr$ for which   $\chi'$  is never zero.  To prove this fact
assume that $f=x_1$ is a local coordinate.  Since $g$ is constant on the leaves 
$\{x_1={\rm constant}\}$,  $g$ must be independent of $x_2,...,x_n$, i.e., $g=\chi(x_1)$.
Since $dg$ is never zero, $\chi'$ is never zero.  Finally, 
$$
dd^\f g\ =\ \chi'(f)dd^\f f +\chi''(f) df\wedge d^\f f
\ =\ dg\wedge\left(\a_f+{{\chi''(f) }\over{\chi'(f) }}d^\f f\right) +\chi'(f) \sigma_f
$$
which proves that if $f$ is \fh mod $d$, then $g=\chi(f)$ is also. \qed

\medskip\Prop{\II.3}  {\sl    If $f$ is \fh mod $d$, then each (non-critical) hypersurface $\{f=C\}$ is \ffl.}

\pf Suppose $\x\in G(\f)$ is tangent to $\{f=C\}$, i.e., $\nabla f \hk \x=0$.  Then
$
(dd^\f f)(\x) \ =\ (df\wedge \a_f)(\x) + \s_f(\x).
$
Since $\s_f$ is of type $\L(\f)^\perp$, we have $\s_f(\x)=0$, and $(df\wedge\a_f)(\x)= \a_f(\nabla f\hk \x)=0$. \qed\medskip

Recall from Proposition \AA.13 that for any $f\in C^\infty$ and any $\f$-submanifold $M$,
we have
$$
(dd^\f f-df\wedge \a_f)\bigr|_M\ =\ *_M(\D_M f) - d(f\bigr|_M)\wedge\a_f \bigr|_M.
\eqno{(\II.2)}
$$
This proves

\Prop{\II.4}  {\sl  If $f$ is \fh mod $d$ and $M$ is a $\f$-submanifold, then
$u\equiv f\bigr|_M$ satisfies the partial differential equation
$$
\D_M u \ =\ *(du\wedge\beta)   \qquad {\rm on\ \ } M
\eqno{(\II.3)}
$$
where $\b=\a_f\bigr|_M$.}

The maximum principle is applicable to solutions to (\II.3).  See for instance [BJS].

\Cor{\II.5}  {\sl  Suppose $(M,\G)$ is a compact $\f$-submanifold with boundary.
Then for each function $f$ which is \fh mod $d$ and each point $x\in M$, one has}
$$
\inf_\G \ \leq f(x) \ \leq \ \sup_\G f
\eqno{(\II.4)}
$$

\Cor{\II.6}  {\sl   Suppose $(M,\G)$ is as above.  
If $\G\subset \{f=C\}$, then $M\subset \{f=C\}$.}

\Prop{\II.7}  {\sl  Suppose $(M,\G)$ is a compact $\f$-submanifold with boundary, and suppose
$f$ is a function on $X$ which is \fh mod $d$. If $f$ is constant on $\G$, then}
$$
d^\f f\bigr|_\G \ \equiv\ 0 \qquad{\rm (pointwise).}
\eqno{(\II.5)}
$$

\pf
By Corollary \II.6, $f$ is constant on $M$. We then apply the following.

\Lemma{\II.8}  {\sl For any function $f$ constant on $M$, $d^\f f\bigr|_\G \ \equiv\ 0$.}

\pf
At $x\in\G$, we have $\overrightarrow M = e\wedge \overrightarrow \G$ for some
$e$ tangent to $M$.  Since $f$ is constant on $M$, $\nabla f\perp \span \overrightarrow M$.
Now $(d^\f f)(\overrightarrow\G)= (\nabla f \hk \f)(e\hk \overrightarrow M) = \f((\nabla f)\wedge (e\hk\overrightarrow M))=0$ since $\nabla f \wedge (e\hk  \overrightarrow M)$ is a first cousin of 
$\overrightarrow M\in G(\f)$. \qedqed

\medskip

Our next objective is to show that, for the large class of  {\sl normal} calibrations, a function
$f$ is \fh mod $d$ if and only if its level sets are \ffl.

Suppose $\f\in \L^pV$ is a calibration on a euclidean vector space $V$.  For each hyperplane
$W\subset V$, $\f\bigr|_W\in\L^p W$ has comass $\leq 1$ and, in fact, $<1$ if and only if
$G(\f\bigr|_W)$ is empty.

\Def{\II.9}  The calibration $\f\in \L^pV$ is {\bf normal} if, for every hyperplane $W\subset V$
$$
\L(\f\bigr|_W)^\perp\ = \L(\f)^\perp\bigr|_W
$$
as subspaces of  $\L^pW$.  A calibration $\f$ on a manifold $X$ is {\bf normal}
if $\f_x\in\L^p T_xX$ is normal for each $x\in X$.

\Prop{\II.10}  {\sl  Suppose $\f$ is a normal calibration on $X$, and $f\in C^\infty(X)$ has a         never-vanishing gradient.  Then
\smallskip

\centerline{\sl  $f$ is \fh mod $d$ }
\smallskip
\noindent
if and only if }

\smallskip

\centerline{\sl  each level set $\{f=C\}$  is \ffl.}
\smallskip

\pf Suppose each  level set $\{f=C\}$  is \ffl.  That is
$$
(dd^\f f)(\x) \ =\ 0 \qquad {\rm for\ all\ }\x\in G(\f) \ \ {\rm which\  are\ tangential\ to\ }
\{f=C\}.
\eqno{(\II.6)}
$$
Note that at a point $x\in X$, $G(\f\bigr|_W) = \{\x\in G(\f) : \x $ is tangential to $W\}$.  Let
$W=\ker df$. Then (\II.6) is equivalent to 
$$
dd^\f f\bigr|_W \in \L\left(\f\bigr|_W\right)^\perp
\eqno{(\II.7)}
$$
Now $f$ is \fh mod $d$ if 
$$
dd^\f f\ =\ df\wedge \a_f +\s_f  \qquad \s_f \in \L(\f)^\perp
\eqno{(\II.8)}
$$
or equivalently
$$
dd^\f f\bigr|_W \in \L(\f)^\perp \bigr |_W
\eqno{(\II.9)}
$$
If $\f$ is normal, then $$ \L\left(\f\bigr|_W\right)^\perp  \subset  \L(\f)^\perp \bigr |_W$$
and (\II.7) implies (\II.9).
\qed

\Prop{\II.11}  {\sl The following calibrations are normal.

\medskip
 
\qquad 1.  A K\"ahler or $p$th power K\"ahler calibration.

\medskip
 
\qquad 2.  A Special Lagrangian calibration.

\medskip
 
\qquad 3.  An associative, coassociative or Cayley calibration.

\medskip
 
\qquad 4.  A quaternionic calibration.

}

\vfill\eject

\vskip .3in



\centerline{\bf References}

\vskip .2in

\noindent
[Al]   S. Alesker,  {\sl  Non-commutative linear algebra and  plurisubharmonic functions  of quaternionic variables}, Bull.  Sci.  Math., {\bf 127} (2003), 1-35. also ArXiv:math.CV/0104209.  

\smallskip

\noindent
[AV]   S. Alesker and M. Verbitsky,  {\sl  Plurisubharmonic functions  on hypercomplex manifolds and HKT-geometry}, arXiv: math.CV/0510140  Oct.2005

\smallskip

\noindent
[Alex]   H. Alexander,  {\sl  Holomorphic chains and the support hypothesis conjecture}, J. Amer. Math.
Soc., {\bf 10} (1997), 123-138.

\smallskip

\noindent
[A] F. J. Almgren, Jr.,  {\sl  $Q$-valued functions minimizing Dirichlet's integral and the regularity of area-minimizing rectifiable currents up to codimension 2},  World Scientific Monograph Series
in Mathematics, 1, World Scientific Publishing Co.River Edge, NJ, 2000.

\smallskip

\noindent
[BJS]   L. Bers, F. John and M. Schechter,  {Partial Differential Equations}, Interscience, J. Wiley,
1964.

\smallskip

\noindent
[Be]   A. Besse,  {Einstein Manifolds}, 
1987.

\smallskip

\noindent
 [B]   R. L.  Bryant,  {\sl
 Calibrated cycles of codimension 3 in compact simple Lie groups},
  to  appear.
 \smallskip

\noindent
 [BS]   R. L.  Bryant and S. M. Salamon,  {\sl
 On the construction of some complete metrics with exceptional holomony},
  Duke Math. J. {\bf 58} (1989),   829-850.

 \smallskip

\noindent
 [C]   E.  Calabi,  {\sl
 M\'etriques k\"ah\'eriennes et fibr\'es holomorphes},
 Annales scientifiques de l'\'Ecole Normale Superieure {\bf 12} (1979),   269-294.

 \smallskip

\noindent
 [Ch]   S.-S. Chern,  {\sl
 On a generalization of K\"ahler geometry},  Lefschetz Jubilee Volume, Princeton Univ. Press,
 (1957), 103-121.

 \smallskip

 \noindent
[DS] J. Duval and N. Sibony, {\sl
Polynomial convexity, rational convexity and currents},
  Duke Math. J. {\bf 79}  (1995),     487-513.

 \smallskip

\noindent
 [C]   T. Eguchi and A. J. Hanson,  {\sl
Asymptotically flat solutions to Euclidean gravity},
Physics Letters  {\bf 74B} (1978),   249-251.

 \smallskip

\noindent
[F]   H. Federer, Geometric Measure  Theory,
 Springer--Verlag, New York, 1969.

 \smallskip

\noindent
 [FF]   H. Federer and W. Fleming, {\sl
Normal and Integral Currents},
Annals of Math.  {\bf 72} (1960),   458-520.

 \smallskip

\noindent
 [Fu]   L. Fu, {\sl  On the boundaries of Special Lagrangian submanifolds},
Duke Math. J.   {\bf 79}   no. 2 (1995),   405-422.

 \smallskip

\noindent
[J]   D. D. Joyce,     
{    Compact Manifolds with Special Holonomy},    
Oxford University Press, Oxford, 2000.

\smallskip

\noindent
 [GL]  K. Galicki and B. Lawson, {\sl  Quaternionic reduction and quaternionic orbifolds}, Math. Ann. {\bf 282} (1989), 1-21.

 \smallskip

\noindent
[GR]   H. Grauer and R. Remmert,  Coherent Analytic Sheaves, Springer-Verlag, Berlin-Heidelberg, 1984.
\smallskip

\noindent
[GZ]  V. Guedj and A. Zeriahi,     
{\sl    Intrinsic capacities on compact K\"ahler manifolds},    
Preprint Univ. de Toulouse , 2003

\smallskip

\noindent
[H$_1$]  F.R. Harvey,
Holomorphic chains and their boundaries, pp. 309-382 in ``Several Complex
Variables, Proc. of Symposia in Pure Mathematics XXX Part 1'', 
A.M.S., Providence, RI, 1977.

\noindent
[H$_2$]  F.R. Harvey,
Spinors and Calibrations,  Perspectives in Mathematics, vol. 9 Academic Press, Boston, 1990

\noindent
[HK] F. R. Harvey and  A.  W. Knapp, {\sl  Positive (p,p)-forms, Wirtinger's inequality and currents}, 
Value-Distribution Theory, Part A (Proc. Tulane Univ. Program on Value-Distribution Theory
in Complex Analysis and Related Topics in Differential Geometry, 1972-73),  pp. 43-62,
Dekker, New York, 1974.

 \smallskip

\noindent
[HL$_1$] F. R. Harvey and H. B. Lawson, Jr, {\sl On boundaries of complex
analytic varieties, I}, Annals of Mathematics {\bf 102} (1975),  223-290.

 \smallskip

\noindent
[HL$_2$] F. R. Harvey and H. B. Lawson, Jr, {\sl On boundaries of complex
analytic varieties, II},  Annals of Mathematics {\bf 106} (1977), 
213-238.

 \smallskip

   \noindent 
 {[HL$_3$]} F. R. Harvey and H. B. Lawson, Jr, {\sl Calibrated geometries},  Acta Mathematica 
{\bf 148} (1982), 47-157.

 \smallskip

 
   \noindent 
 {[HL$_4$]} F. R. Harvey and H. B. Lawson, Jr, {\sl Boundaries of positive holomorphic chains in K\"ahler manifolds}, Stony Brook Prerprint, 2005.

 \smallskip

\noindent
[HLZ] F. R. Harvey, H. B. Lawson, Jr. and J. Zweck, {\sl A
deRham-Federer theory of differential characters and character duality},
Amer. J. of Math.  {\bf 125} (2003), 791-847. ArXiv:math.DG/0512251

 \smallskip

\noindent
[HP] F. R. Harvey, J. Polking, {\sl Extending analytic objects},
Comm.  Pure Appl. Math. {\bf 28} (1975), 701-727.

 \smallskip

\noindent
[HW$_1$] F. R. Harvey,  R. O. Wells, Jr.,  {\sl Holomorphic approximation and hyperfunction 
theory on a $C^1$ totally real submanifold of a complex manifold},
  Math.  Ann. {\bf 197} (1972),  287-318.

 \smallskip

\noindent
[HW$_2$] F. R. Harvey,  R. O. Wells, Jr.,  {\sl Zero sets of non-negatively strictly plurisubharmonic
functions},
  Math.  Ann. {\bf 201} (1973),  165-170.

 \smallskip

   \noindent
[HS]    F.R. Harvey and B. Shiffman,    {\sl  A characterization of
holomorphic chains},    Ann. of Math.,
 {\bf 99}  (1974), 553-587.

\smallskip

   \noindent
[K]    J. King,    {\sl   The currents defined by analytic varieties},  Acta Math.   {\bf 127}  no. 3-4 (1971),
185-220.

\smallskip

   \noindent
[L$_1$]    H. Blaine Lawson, Jr.,      Minimal Varieties in Real and Complex Geometry, Les Presses de 
L'Universite de Montreal,  1974.

 \smallskip

   \noindent
[L$_2$]    H. Blaine Lawson, Jr.,    {\sl  Minimal Varieties},    
Proceedings of Symposia in Pure Mathematics 
{\bf 27} (1974), 61-93.

 \smallskip

   \noindent
[L$_3$]    H. Blaine Lawson, Jr.,    {\sl  The stable homology of a flat torus},    
Mathematica Scandinavica 
{\bf 36} (1975), 49-73.

\smallskip

 \noindent 
 {[L$_4$]} {\sl Geometric aspects of the genealized Plateau problem}, 
pp. 7-13 in  Proceedings  of the International Congress of Mathematicians, 1974,
vol.2, Vancouver, Canada, 1975.

 \smallskip

 \noindent 
 {[L$_5$]} {\sl The question of holomorphic carriers}, Proceedings of
Symposia in  Pure Mathematics {\bf 30}, American Mathematical 
Society (1976), 115-124.

 \smallskip

   \noindent
[S]   H. H. Schaefer,  Topological Vector Spaces,    Springer Verlag,
New York,  1999.

\smallskip

   \noindent
[Sh]      B. Shiffman,    {\sl Complete characterization of holomorphic chains of codimension one},    
Math. Ann. {\bf 274} (1986), 233-256.
\smallskip

\noindent
[U]  I.  Unal, Ph.D. Thesis, Stony Brook, 2006.

   \noindent
[V]    M. Verbitsky.,    {\sl  Manifolds with parallel differential forms and K\"ahler identities
for $G_2$-manifolds},  arXiv : math.DG/0502540 (2005).
 \smallskip

\end